\documentclass[11pt]{article}
\usepackage{amsfonts}
\usepackage{amsfonts,latexsym,amsmath,amscd,geometry}
\usepackage{amssymb}
\usepackage{latexsym}

\newcommand \nc{\newcommand}
\newtheorem{theorem}{Theorem}[section]
\newtheorem{lemma}[theorem]{Lemma}

\newtheorem{corollary}[theorem]{Corollary}
\newtheorem{definition}[theorem]{Definition}

\newtheorem{remark}[theorem]{Remark}

\nc{\ba}{\begin{array}}\nc{\ea}{\end{array}}
\nc{\be}{\begin{eqnarray}}\nc{\ee}{\end{eqnarray}}
\nc{\beq}{\begin{equation}}\nc{\eeq}{\end{equation}}
\nc{\bex}{\begin{eqnarray*}}\nc{\eex}{\end{eqnarray*}}
\nc{\btm}{\begin{theorem}} \nc{\etm}{\end{theorem}}
\nc{\blm}{\begin{lemma}} \nc{\elm}{\end{lemma}}
\nc{\R}{\mathbb{R}} \nc{\va}{\varepsilon} \nc{\ls}{\limits}

\def\de{\Delta}\def\pa{\partial}\def\om{\Omega}
\def\pf{\noindent{\bf Proof.\quad}}\def\endpf{\hfill$\Box$}
\def\les{\lesssim}\def\u{\dot{u}}\def\di{\mbox{div\,}}

\newcommand \qed {\hfill $\Box$}

\begin{document}
\title{Strong solutions of the compressible nematic liquid crystal flow}
\author{Tao Huang\footnote{Department of Mathematics, University of Kentucky,
Lexington, KY 40506, USA. } \quad Changyou Wang$^*$ \quad Huanyao Wen\footnote{School of Mathematical Sciences, South China Normal University, Guangzhou, 510631, P. R. China.}}

\maketitle

\begin{abstract} We study strong solutions of the simplified Ericksen-Leslie system modeling  compressible nematic liquid crystal flows in a domain $\Omega\subset\mathbb R^3$. We first prove the local existence
of unique strong solutions provided that the initial data $\rho_0, u_0, d_0$ are sufficiently regular and satisfy a natural compatibility condition. The initial densitiy  function $\rho_0$ may vanish on an open subset (i.e., an initial
vacuum may exist). We then prove a criterion for possible breakdown of such a local strong solution at finite
time in terms of blow up 
of the quantities $\displaystyle\|\rho\|_{L^\infty_tL^\infty_x}$ and $\displaystyle\|\nabla d\|_{L^3_tL^\infty_x}$. 
\end{abstract}

\section {Introduction}
\setcounter{equation}{0}
\setcounter{theorem}{0}

Nematic liquid crystals are aggregates of molecules which possess same orientational order and are made of elongated,
rod-like molecules. The continuum theory of liquid crystals was developed by Ericksen \cite{Ericksen}  and Leslie \cite{Leslie}
during the period of 1958 through 1968, see also the book by de Gennes \cite{Gennes}.  Since then there have been remarkable
research developments in liquid crystals from both theoretical and applied aspects.  When the fluid containing nematic liquid crystal materials is at rest, we have the well-known Ossen-Frank theory for static nematic liquid crystals, see Hardt-Lin-Kinderlehrer \cite{HLK} on the analysis of energy minimal configurations of namatic liquid crystals.
 In general, the motion of fluid always takes place. The so-called Ericksen-Leslie system is a macroscopic continuum description of the time evolution of the materials under the influence of both the flow velocity field $u$ and the macroscopic description of the microscopic orientation configurations $d$ of rod-like liquid crystals. 

When the fluid is an incompressible, viscous  fluid, Lin \cite{Lin} first derived a simplified Ericksen-Leslie equation modeling liquid crystal flows in 1989. Subsequently, Lin and Liu \cite{LL, LL1} made some important analytic studies, such as the existence of weak and strong solutions and the partial regularity of suitable solutions,  of the simplified Ericksen-Leslie system, under the assumption that the liquid crystal director field is of varying length by Leslie's terminology or variable degree of orientation by Ericksen's terminology.  

When the fluid is allowed to be compressible, the Ericksen-Leslie system becomes more complicate and there seems very few analytic works available yet.  We would like to mention that very recently, there have been both modeling study, see Morro \cite {Morro},  and numerical study, see Zakharov-Vakulenko \cite{ZV}, on the hydrodynamics of compressible nematic liquid crystals  under the influence of temperature gradient or electromagnetic forces.   

This paper, and the companion paper \cite{huang-wang-wen}, aims to study the strong solutions of the flow of compressible nematic liquid crystals and the blow up criterions.

Let $\Omega\subset\mathbb R^3$ be a domain. We will consider the simplified version of Ericksen-Leslie system modeling the flow of compressible nematic liquid crystals in $\Omega$:
\begin{align}
\rho_t+\nabla \cdot (\rho u)&=0, \label{clc-1} \\
\rho u_t+\rho u\cdot\nabla u+\nabla (P(\rho))&=\mathcal L u-\nabla d\cdot \de d,\label{clc-2}\\
d_t+u\cdot\nabla d&=\de d+|\nabla d|^2d,\label{clc-3}
\end{align}
where $\rho:\om\times[0,+\infty)\to\mathbb
R^1$ is the density function of the fluid, \ $u:\om\times[0,+\infty)\to\mathbb R^3$ represents
velocity field of the fluid,  $P=P(\rho)$
represents the pressure function, 
$d:\om\times[0,+\infty)\to S^2$ represents the macroscopic
average of the nematic liquid crystal orientation field, $\nabla\cdot$ is the divergence operator in $\mathbb R^3$,
and $\mathcal L$ denotes the Lam$\acute{\rm e}$ operator:
$$\mathcal Lu=\mu\de u+(\mu+\lambda)\nabla\mbox{div }u,$$
where $\mu$ and $\lambda$ are shear viscosity and the bulk viscosity coefficients 
of the fluid respectively that satisfy the physical condition:
\beq\label{viscosity} \mu>0,\quad  2\mu+3\lambda\geq 0. \eeq
We refer to the readers to consult the recent preprint \cite{ding-huang-wen-zi} by Ding-Huang-Wen-Zi 
for the derivation for the system (\ref{clc-1})-(\ref{clc-3}) based on energetic-variational approaches.
Throughout this paper, we assume that
\beq\label{regularity_p}
P:[0,+\infty)\to\mathbb R \ {\rm{is\ a\ locally\ Lipschitz\ continuous\ function}}.
\eeq
Notice that (\ref{clc-1}) is the equation of conservation of mass, (\ref{clc-2}) is the equation of linear momentum, and
(\ref{clc-3}) is the equation of angular momentum. 
We would like to point out that the system (\ref{clc-1})-(\ref{clc-3})
includes several important equations as special cases:

(i) When $\rho$ is constant, the equation (\ref{clc-1}) reduces
to the incompressibility condition of the fluid ($\nabla\cdot u=0$), and the system (\ref{clc-1})-(\ref{clc-3}) becomes
the equation of incompressible flow of namatic liquid crystals provided that $P$ is a unknown pressure function. This
was previously proposed by Lin \cite{Lin} as a simplified Ericksen-Leslie equation modeling incompressible liquid crystal flows.

(ii) When $d$ is a constant vector field,  the system (\ref{clc-1})-(\ref{clc-2})
becomes a compressible Navier-Stokes equation, which is an extremely important equation to describe 
compressible fluids (e.g., gas dynamics).
It has attracted great interests among many analysts and there have been many important developments (see,
for example,  Lions \cite{Lions2}, Feireisl \cite{Feireisl} and references therein).

(iii) When both $\rho$ and $d$ are constants, the system (\ref{clc-1})-(\ref{clc-2}) becomes the incompressible Naiver-Stokes
equation provided that $P$ is a unknown pressure function, the fundamental equation to describe Newtonian fluids (see, for example, Lions \cite {Lions1} and Temam \cite{Temam} for  survey of important developments).

(iv) When $\rho$ is constant and $u=0$, the system (\ref{clc-1})-(\ref{clc-3}) reduces to the equation for heat flow of harmonic maps into $S^2$.  There have been extensive studies on the heat flow of harmonic maps in the past few decades (see, for example,
the monograph by Lin-Wang \cite{LW1} and references therein).

From the viewpoint of partial differential equations, the system (\ref{clc-1})-(\ref{clc-3}) is a highly nonlinear system coupling between hyperbolic equations and parabolic equations.
It is very challenging to understand and analyze such a system, especially when the density function $\rho$ may vanish
or the fluid takes vacuum states.

In this paper, we will consider the following initial condition:
\begin{align}
(\rho, u, d)\Big|_{t=0}&= (\rho_0, u_0, d_0), \label{clcinitial}
\end{align}
and one of the three types of boundary conditions:

(1) Cauchy problem:
\begin{align}\label{clcboundary1}
\om=\R^3, \ \mbox{ and }\rho,\ u \mbox{ vanish at infinity and\ } d \mbox{ is\ constant  at infinity (in some weak
sense)}.
\end{align}

(2) Dirichlet and Neumann boundary condition for $(u,d)$: $\om\subset\R^3$ is a bounded smooth domain, and
\begin{align}
(u,\ \frac{\partial d}{\partial\nu})\Big|_{\pa\Omega}&= 0, \label{clcboundary2}
\end{align}
where $\nu$ is the unit outer normal vector of $\pa\om$.

(3) Navier-slip and Neumann boundary condition for $(u,d)$: $\om\subset\R^3$ is bounded, simply connected, smooth domain, and
\begin{align}
(u\cdot\nu,(\nabla\times u)\times\nu,\ \frac{\partial d}{\partial\nu})\Big|_{\pa\Omega}&= 0, \label{clcboundary3}
\end{align}
where $\nabla\times u$ denotes the vorticity field of the fluid.

To state the definition of strong solutions to the initial and boundary value
problem  (\ref{clc-1})-(\ref{clc-3}), (\ref{clcinitial}) together with (\ref{clcboundary1})
or (\ref{clcboundary2}) or (\ref{clcboundary3}), we introduce some notations. 

We denote 
$$\int f\,dx =\int_\Omega f \,dx.$$
For $1\le r\le \infty$, denote the $L^r$ spaces and the standard Sobolev spaces as follows:
$$L^r=L^r(\Omega),  \ D^{k,r}=\left\{ u\in L^1_{\rm{loc}}(\Omega): \|\nabla^k u \|_{L^r}<\infty\right\},$$
$$W^{k,r}=L^r\cap D^{k,r},  \ H^k=W^{k,2}, \ D^k=D^{k,2},$$
$$D_0^1=\Big\{u\in L^6: \ \|\nabla u\|_{L^2}<\infty,\
{\rm{ and \ satisfies}}\ (\ref{clcboundary1}) \ {\rm{or}}\ (\ref{clcboundary2}) \ {\rm{or}}\ (\ref{clcboundary3}) 
\ {\rm{ for\ the\ part\ of }}\ u\Big\},$$
$$H_0^1=L^2\cap D_0^1, \ \|u\|_{D^{k,r}}=\|\nabla^k u\|_{L^r}.$$
Denote $$Q_T=\Omega\times [0,T] \ (T>0),$$
and let
$$\mathcal D(u)=\frac12\left(\nabla u+(\nabla u)^t\right)$$
denote the deformation tensor, which is the symmetric part of the velocity gradient.

\begin{definition} For $T>0$, $(\rho, u,d)$ is called a strong solution to the compressible nematic liquid crystal flow (\ref{clc-1})-(\ref{clc-3}) in $\Omega\times (0,T]$, if for some $q\in (3, 6]$, 
\bex &0\le \rho\in
C([0,T];W^{1,q}\bigcap H^1),\ \rho_t\in C([0,T];L^2\bigcap L^q);&\\
& u\in C([0,T];D^2\bigcap D^1_0)\bigcap L^2(0,T;D^{2,q}),\ u_t\in
L^2(0,T;D^1_0),\
 \sqrt{\rho} u_t\in
L^\infty(0,T;L^2);&\\
& \nabla d\in C([0,T];H^2)\bigcap L^2(0,T;
H^3),\ d_t\in C([0,T];H^1)\bigcap L^2(0,T; H^2),\
 |d|=1 \ {\rm{in}}\ \overline{Q}_T;& \eex
and $(\rho,u, d)$ satisfies (\ref{clc-1})-(\ref{clc-3}) a.e. in $\Omega\times (0,T]$.

\end{definition}

The first main result is concerned with local existence of strong
solutions.

\begin{theorem}\label{App-thm:local}
Assume that $P$ satisfies (\ref{regularity_p}), $\rho_0\geq0$, $\rho_0\in W^{1,q}\bigcap H^1\bigcap L^1$ for some
$q\in(3,6]$, $u_0\in D^2\bigcap D^1_0$, $\nabla d_0\in H^2$ and
$|d_0|=1$ in $\overline{\Omega}$. If, in additions, the following compatibility
condition 
\be\label{first3.1} \mathcal Lu_0-\nabla
(P(\rho_0))-\de d_0\cdot\nabla d_0=\sqrt{\rho_0}g\
{\rm{ for\ some}}\ g\in L^2(\Omega, \mathbb R^3)
\ee
holds,  then there exist a positive time $T_0>0$
 and a unique strong solution $(\rho,u, d) $ of (\ref{clc-1})-(\ref{clc-3}),
(\ref{clcinitial}) together with (\ref{clcboundary1}) or 
(\ref{clcboundary2}) or (\ref{clcboundary3}) in $\Omega\times (0,T_0]$.

\end{theorem}

We would like to point out that  an analogous existence theorem of local strong solutions to the isentropic\footnote{i.e. $P=a\rho^\gamma$ for some $a>0$ and $\gamma>1$.} compressible Naiver-Stokes equation, under the first two boundary conditions (\ref{clcboundary1}) and
(\ref{clcboundary2}), has been previously established  by Choe-Kim \cite{CK} and  Cho-Choe-Kim \cite{CCK}. A byproduct of our theorem \ref{App-thm:local} also yields the existence of local strong solutions to a larger class of compressible Navier-Stokes
equations under the Navier-slip boundary condition (\ref{clcboundary3}), which seems not available in the literature.  

In dimension one, Ding-Lin-Wang-Wen \cite{DLWW} have proven that the local strong solution to
(\ref{clc-1})-(\ref{clc-3}) under (\ref{clcinitial}) and (\ref{clcboundary2}) is global.  For dimensions at least two, it is reasonable to believe that the local strong solution to (\ref{clc-1})-(\ref{clc-3}) may cease to exist globally.
In fact, there exist finite time singularities of the (transported) heat flow of harmonic maps (\ref{clc-3}) in
dimensions two or higher (we refer the interested readers to \cite{LW1} for the exact references). An important question to ask would be what is the main mechanism of possible break down of local strong (or smooth)
solutions. 

Such a question has been studied for the incompressible Euler equation or the Navier-Stokes equation by Beale-Kato-Majda in their poineering work \cite{BKM}, which showed that the $L^\infty$-bound of vorticity $\nabla\times u$ must blow up. Later, Ponce \cite{Ponce} rephrased the BKM-criterion in terms of the deformation tensor $\mathcal D(u)$.  

When dealing with the isentropic compressible Navier-Stokes equation, there have recently been several very interesting works on the blow up criterion. For example,
if $0<T_*<+\infty$ is the maximum time for strong solution, then (i) Huang-Li-Xin \cite{HLX0}  established
a Serrin type criterion: 
$\lim_{T\uparrow T_*} \big(\|{\rm{div}} u\|_{L^1(0,T; L^\infty)}+\|\sqrt{\rho}  u\|_{L^s(0,T; L^r)}\big)=\infty$ for $\frac{2}{s}+\frac{3}{r}\le 1, \ 3<r\le\infty$; (ii) Sun-Wang-Zhang \cite{Sun-Wang-Zhang}, and independently \cite{HLX0}, showed that if $7\mu>\lambda$, then 
$\lim_{T\uparrow T_*} \|\rho\|_{L^\infty(0,T; L^\infty)}=\infty$; and (iii) Huang-Li-Xin \cite{HLX} showed
that $\lim_{T\uparrow T_*} \|\mathcal D(u)\|_{L^1(0,T; L^\infty)}=\infty$.  
 
When dealing the heat flow of harmonic maps (\ref{clc-3}) (with $u=0$),  Wang \cite{Wang} obtained
a Serrin type regularity theorem, which implies that if $0<T_*<+\infty$ is the first singular time for
local smooth solutions, then $\lim_{T\uparrow T_*} \|\nabla d\|_{L^2(0,T; L^\infty)}=\infty$.

When dealing with the incompressible nematic liquid crystal flow,  Lin-Lin-Wang \cite{Lin-Lin-Wang}  and Lin-Wang \cite{LW} have established the global existence of a unique "almost strong" solution\footnote{that has at most finitely many possible singular time.} for the initial-boundary value problem in  bounded domains
in dimension two, see also Hong \cite{Hong} and Xu-Zhang \cite{Xu-Zhang} for some related works.   
In dimension three,  for the incompressible nematic liquid crystal flow Huang-Wang \cite{HW} have obtained a BKM type blow-up criterion very recently,  while the existence of global weak solutions still
remains to be a largely open question. 

Motivated by these works on the blow up criterion of local strong solutions to the Navier-Stokes equation and
the incompressible nematic liquid crystal flow, we will establish in this paper the following blow-up criterion 
of breakdown of local strong solutions under the boundary condition (\ref{clc-1}) or (\ref{clc-2}). 

\btm\label{umaintheorem}{\it 
Let $(\rho, u, d)$ be a strong solution of the initial boundary problem
(\ref{clc-1})-(\ref{clc-3}), (\ref{clcinitial}) together with (\ref{clcboundary1})
or (\ref{clcboundary2}). Assume that $P$ satisfies (\ref{regularity_p}),
and  the initial data $(\rho_0, u_0,d_0)$
satisfies (\ref{first3.1}). If $0<T_*<+\infty$ is the maximum time of existence
and $7\mu>9\lambda$,
then  
\beq\label{clcblpcondition}
\lim_{T\uparrow T_*} \Big(\|\rho\|_{L^{\infty}(0,T; L^\infty)}+\|\nabla
d\|_{L^3(0,T;L^{\infty})}\Big)=\infty. \eeq}\etm

We would like to make a few comments of Theorem \ref{umaintheorem}.

\begin{remark} {\rm{(a) Since we can't yet prove Lemma 4.2  for the Navier-slip and Neumann boundary condition (\ref{clcboundary3}),  it is unclear whether Theorem \ref{umaintheorem} remains to be true
under the boundary condition (\ref{clcboundary3}). \\
(b) In \cite{huang-wang-wen}, we obtained a blow-up criterion of (\ref{clc-1})-(\ref{clc-3}) under the initial condition (\ref{clcinitial}) and the boundary condition (\ref{clcboundary1})
or (\ref{clcboundary2}) or (\ref{clcboundary3}) in terms of $u$ and $\nabla d$: if
$0<T_*<+\infty$ is the maximum time of existence of strong solutions, then
$$\lim_{T\uparrow T_*} \Big(\|\mathcal D(u)\|_{L^1(0,T; L^\infty)}+\|\nabla d\|_{L^2(0,T; L^\infty)}\Big)=+\infty.$$
(b) For compressible liquid crystal flows without the nematicity constraint ($|d|=1$)\footnote{the right hand side of equation (\ref{clc-3}) is replaced by $\Delta d+f(d)$ for some smooth function $f:\mathbb R^3
\to\mathbb R^3$, e.g. $f(d)=(|d|^2-1)d$.}, Liu-Liu \cite{Liu-Liu}
have recently obtained a Serrin type criterion on the blow-up of strong solutions.\\
(c) It is a very interesting question to ask whether there exists a global weak solution to the initial-boundary value problem of (\ref{clc-1})-(\ref{clc-3}) in dimensions at least two. In dimension one, such an
existence has been obtained by Ding-Wang-Wen \cite{DWW}.}}
\end{remark}

Now we briefly outline the main ideas of the proof, some of which are inspired by earlier works on the isentropic
compressible Navier-Stokes equations by \cite{CCK}, \cite{Sun-Wang-Zhang}, and \cite{HLX}. To obtain the existence of a unique local strong solution to (\ref{clc-1}-(\ref{clc-3}), under (\ref{clcinitial}) and
(\ref{clcboundary1}) or (\ref{clcboundary2}) or (\ref{clcboundary3}), we employ the Galerkin's method
that requires us to establish a {\it priori}  estimate of the quantity 
$$\|\rho(t)\|_{H^1\cap W^{1,q}}+\|\nabla u(t)\|_{L^2}+\|\sqrt{\rho} u_t(t)\|_{L^2}
+\|\nabla^2 d(t)\|_{L^2}, \ 3<q\le 6 $$
for strong solutions $(\rho, u,d)$  in the form of a  Gronwall type inequality. See Theorem 2.1.
It may be of independent interest that we establish $W^{2,q}$-estimate for the Lam\'e equation under the Navier-slip boundary condition, see Lemma 3.1.

To prove  the blow-up criterion (\ref{clcblpcondition}) of Theorem \ref{umaintheorem} in terms of $\rho$ and $\nabla d$, a critical step is to establish the $\displaystyle L^\infty_tL^q_x$-estimate of $\nabla\rho$. From the continuity equation (\ref{clc-1}), this requires that the Lipschitz norm of velocity field $u$, or 
$\displaystyle\|\nabla^2 u(t)\|_{L^q}$ 
is bounded in $L^1_t$. This is done in several steps. \\
(1) We show that under the condition $7\mu>9\lambda$,
the bound of $\displaystyle(\|\rho\|_{L^\infty_tL^\infty_x}+\|\nabla d\|_{L^3_tL^\infty_x})$ in equations
(\ref{clc-2}) and (\ref{clc-3}) can yield both a high integrability and a high order estimate of $u$ and $\nabla d$, i.e. both  $\displaystyle (\|\rho^\frac15 u\|_{L^\infty_tL^5_x}+\|\nabla d\|_{L^\infty_tL^5_x})$ 
and $\displaystyle(\|\nabla u\|_{L^\infty_tL^2_x}+\|\nabla^2 d\|_{L^\infty_tL^2_x})$ are bounded.
See Lemma 4.2. \\
(2) Based on these estimates from (1), we establish that $\nabla^3 d$ is bounded in $L^\infty_tL^2_x$ and $\nabla u$ is bounded in $\displaystyle L^2_tW^{1,q}_x+L^\infty_t ({\rm{BMO}}_x)$. To achieve it, we adapt  the approach, due to Sun-Wang-Zhang \cite{Sun-Wang-Zhang}, by decomposing $u=w+v$, where $v\in H^1_0(\Omega)$ solves the Lam\'e equation $\mathcal Lv=\nabla (P(\rho))$.
One can prove that $\nabla v\in L^\infty_t ({\rm{BMO}}_x)$ by the elliptic regularity theory. 
The difficult part is  to show that $\nabla^2 w\in L^2_tL^q_x$ for $3<q\le 6$. In order to obtain this 
estimate, we first establish that $\displaystyle(\|\sqrt{\rho}\dot{u}\|_{L^\infty_tL^2_x}+\|\nabla d_t\|_{L^\infty_tL^2_x})$ and $\displaystyle(\|\nabla \dot u\|_{L^2_tL^2_x}+\|d_{tt}\|_{L^2_tL^2_x})$
are bounded by viewing (\ref{clc-2}) as an evolution equation of the material derivative $\dot{u}\equiv u_t+u\cdot\nabla u$ and performing second order energy estimates of both equations (\ref{clc-2}) and (\ref{clc-3}).  Then we employ $W^{2,q}$-estimate of the Lam\'e equation to control $\displaystyle\|\nabla^2 w\|_{L^q}$. The details are illustrated by Lemma 4.4 and Corollary 4.5.\\
(3) We show that $\displaystyle \|\nabla \rho\|_{L^2\cap L^q}$ is bounded by an argument similar to \cite{Sun-Wang-Zhang} \S5.  Then we apply $W^{2,q}$-estimate of the Lam\'e equation again to control
$\displaystyle {\|\nabla^2 u\|_{L^\infty_tL^2_x}}$ and $\displaystyle{ \|u\|_{L^\infty_tD^{2,q}_x}}$. See Lemma 4.6, Corollary 4.7, and Corollary 4.8. 

It is interesting to notice that during the proof of both the existence of a unique local strong solutions and the blow-up criterion for strong solutions, specific forms of the pressure function $P(\rho)$ play no roles and it is the local Lipschitz regularity of $P$ that matters.

The paper is written as follows. In \S2, we derive some a priori estimates for strong solutions or approximate solutions via the Galerkin's method. In \S3, we prove both the local existence by the Gakerlin's method and uniqueness  of strong solutions. In \S4, we discuss the blow up criterion of strong solutions and prove Theorem \ref{umaintheorem}.\\

\noindent{\bf Acknowledgement}. The first two authors are partially supported by NSF grant 1000115.
The work is completed during the visit of third author to University of Kentucky, which is partially supported
by the second author's NSF grant 0600162. The third author wishes to thank the department of Mathematics for its
hospitality.

\section{A {\em priori} estimates}
\setcounter{equation}{0} \setcounter{theorem}{0}

In the section, we will derive some a {\em priori} estimates for strong or smooth
solutions ($\rho$, $u$, $d$) to (\ref{clc-1})-(\ref{clc-3})  on a bounded domain,
associated with the initial condition (\ref{clcinitial}) and the boundary condition (\ref{clcboundary2})
or (\ref{clcboundary3}), provided that the initial density function  has a positive lower bound,
$\rho_0\ge\delta>0$. All these a priori estimates we will obtain are independent of $\delta>0$ and
the size of the domain when $\Omega=B_R$ ($R\ge 1$) is a ball in $\mathbb R^3$, which are the
crucial ingredients to prove the local existence of strong solutions to (\ref{clc-1})-(\ref{clc-3})
when we allow the initial data $\rho_0\ge 0$ and unbounded domain $\Omega=\mathbb R^3$.  Although these estimates may
have their own interests, we mainly apply them to the approximate solutions to (\ref{clc-1})-(\ref{clc-3})
that are constructed by the Galerkin's method. 

Throughout the paper, we denote by $C$  generic constants that depend on
$\|\rho_0\|_{W^{1,q}\cap H^1\cap L^1}$, $\|u_0\|_{D^2\cap D_0^1}$,
$\|\nabla d_0\|_{H^2}$, and $P$, but are independent of $\delta>0$, the solutions $(\rho, u, d)$
and the size of domain when $\Omega=B_R$ ($R\ge 1$) is a ball in $\mathbb R^3$. 
We will also use the obvious notation
$$\|\cdot\|_{X_1\cap\cdots\cap
X_k}=\sum_{i=1}^k\|\cdot\|_{X_i}$$
for Banach spaces $X_i$, $1\le i\le k$ and $k=2,3$. We will use $A\lesssim B$ to denote $A\le CB$ for
some constant generic $C> 0$. 

Let  $(\rho, u, d)$ be a strong solution of (\ref{clc-1})-(\ref{clc-3}) in $\Omega\times (0,T]$
(or the approximate solutions $(\rho^m, u^m, d^m)$ of (\ref{clc-1})-(\ref{clc-3}) constructed by the Galerkin's method in
\S3.2 below). For simplicity, we assume $0<T\le 1$. For  $0<t<T$, set
\begin{equation} \label{phi}
\Phi(t):=\sup\limits_{0\le s\le t}\Big(\|\rho(s)\|_{H^1\cap W^{1,q}}+\|\nabla
u(s)\|_{L^2}+\|\sqrt{\rho}u_t(s)\|_{L^2}+\|\nabla^2
d(s)\|_{H^1}+1\Big).
\end{equation}
The main aim  of this section is to estimate each term of $\Phi$ in terms of some integrals of $\Phi$.
In \S3 below, we will apply arguments of Gronwall's type to prove that $\Phi$ is locally bounded.  

Throughout this section and \S 3, we will let $\mathcal F$ to denote the set that consists of monotonic increasing, locally bounded functions $M$ from
$[0,+\infty)$ to $[0,+\infty)$ with $M(0)=0$, which are independent of $\delta$ and the size of $\Omega$. The reader will see that the exact
form of $M\in\mathcal F$ is not important and may vary from lines to lines during the proof
of the Lemmas. 

Now we state the main theorem of this section.
\begin{theorem}\label{le:Phi} There exists $M\in\mathcal F$ such that  for any $0<t<T$, it holds
 \beq\label{es for Phi}
\begin{split}
\Phi(t)\leq \exp\Big[C\mathcal M(\rho_0, u_0, d_0)+
C\int_0^tM(\Phi(s))\,ds\Big],
\end{split}
\eeq
where 
\begin{equation}\label{rho-u-d}
\mathcal M(\rho_0,u_0,d_0)=1+\left\|\displaystyle\frac{\mathcal Lu_0-\nabla
(P(\rho_0))-\de d_0\cdot\nabla d_0}{\sqrt{\rho_0}}\right\|_{L^2}.
\end{equation}
\end{theorem}

The proof of Theorem \ref{le:Phi} is based on several Lemmas. We may assume $P(0)=0$.  Observe that
 (\ref{regularity_p}) implies that the Lipschitz norm
\beq\label{lip_norm}
B_P(R):=\|P'\|_{L^\infty([0,R])}: [0,+\infty)\to [0,+\infty) \ {\rm{is\ montonic\ increasing\ and \ locally\ bounded}}.
\eeq

\begin{lemma}\label{le:2.1}(energy inequality) There exists $M\in\mathcal F$ such that  for any $0<t<T$, it holds
\beq\label{Energy-identity}\begin{split}&\int_\Omega\left(\rho|u|^2+|\nabla
d|^2\right)\,dx+\int_0^t\int_\Omega\Big[|\nabla u|^2+\left|\Delta
d+|\nabla d|^2d\right|^2\Big]\,dx\le C+\int_0^t M(\Phi(s))\,ds.\end{split} \eeq
\end{lemma}
\pf  Here we only sketch the proof for the boundary condition (\ref{clcboundary3}). Multiplying (\ref{clc-2}) by $u$ and integrating over $\Omega$, using $\displaystyle \de u=\nabla\di
u-\nabla\times(\nabla\times u)$ and (\ref{clc-1}), and applying integration by parts several times,
we obtain
\beq\label{clc2.3}
\frac{1}{2}\frac{d}{dt}\int \rho|u|^2\,dx
+\int(\mu|\nabla\times u|^2+(2\mu+\lambda)|\mbox{div}u|^2)\,dx=\int P(\rho){\rm{div}}u\,dx-\int u\cdot\nabla d\cdot\de d\,dx.
\eeq
Since $\Omega$ is assumed to be simply connected for the boundary condition (\ref{clcboundary3}), we have (see \cite{Wolf}):
\beq\label{div-curl}
\|\nabla u\|_{L^2}\lesssim \|\nabla\times u\|_{L^2}+\|{\rm{div}} u\|_{L^2}, \ \forall u\in H^1(\Omega)
\ {\rm{with}}\ u\cdot\nu=0 \ {\rm{on}}\ \partial\Omega.
\eeq
This and (\ref{viscosity}) imply
\beq
\int (\mu|\nabla\times u|^2+(2\mu+\lambda)|{\rm{div}} u|^2)\,dx\ge \frac{\mu}{3}\int (|\nabla \times u|^2
+|{\rm{div}} u|^2)\,dx\ge \frac{1}{C}\int |\nabla u|^2\,dx.\label{div-curl1}
\eeq
By Cauchy inequality, we have
\beq\label{p-divu}
\left|\int P(\rho) {\rm{div}} u\,dx\right|
\le \frac1{2C}\int |\nabla u|^2\,dx +C\int |P(\rho)|^2\,dx.
\eeq
Multiplying (\ref{clc-3}) by $\de d+|\nabla d|^2 d$ and integrating over $\Omega$, using integration by parts and the fact that
$|d|=1$ we obtain 
\beq\label{clc2.4}
\frac{1}{2}\frac{d}{dt}\int |\nabla d|^2\,dx+\int\left|\de d+|\nabla d|^2d\right|^2\,dx=\int u\cdot\nabla d\cdot\de d\,dx.
\eeq
Combining (\ref{clc2.3}), (\ref{div-curl1}), (\ref{p-divu}),  and (\ref{clc2.4}) together, we obtain
\beq \label{energy_ineq11}
\frac{d}{dt}\int (\rho |u|^2+|\nabla d|^2)\,dx+\int (\frac{1}{C}|\nabla u|^2+|\Delta d+|\nabla d|^2 d|^2)\,dx
\le C\int |P(\rho)|^2\,dx.
\eeq
To estimate the right hand side of (\ref{energy_ineq11}), first observe that by (\ref{lip_norm}) we have\footnote{when $\Omega=B_R$ for $R\ge 1$,
one can the independence of $C$  with respect to $R$ as follows: 
$$\|\rho\|_{L^\infty(B_R)}\leq \max_{x\in B_R} \|\rho\|_{L^\infty(B_1(x))}
\le C\max_{x\in B_R} \|\rho\|_{W^{1,q}(B_1(x))}\le C\|\rho\|_{W^{1,q}(B_R)}.$$}
\beq\label{infty_bound_rho}
\|\rho\|_{L^\infty}+\|P(\rho)\|_{L^\infty}+\|P(\rho)\|_{H^1\cap W^{1,q}}
\le C\Phi+CB_P(\|\rho\|_{L^\infty})\Phi\le M(\Phi)
\eeq
for some $M\in \mathcal F$.
It follows from (\ref{clc-1}) and Sobolev's inequality that
\begin{eqnarray}\label{L2-est-P}
\int |P(\rho)|^2\,dx &=&\int |P(\rho_0)|^2\,dx+2\int_0^t \int P(\rho) P'(\rho)(-\rho {\rm{div}} u -\nabla\rho\cdot u)\,dx\,dt\nonumber\\
&\le& C+C\int_0^t B_P(\|\rho\|_{L^\infty}) (\|P(\rho)\|_{L^3}\|\nabla\rho\|_{L^2}+\|P(\rho)\|_{L^2}\|\rho\|_{L^\infty})\|\nabla u\|_{L^2}\,ds\nonumber\\
&\le& C+\int_0^t M(\Phi(s))\,ds\le C+M(\Phi(t))
\end{eqnarray}
as $M(\Phi(s))$ is increasing and $t\le 1$.
Substituting (\ref{L2-est-P}) into (\ref{energy_ineq11}) and integrating over $[0,t]$ yields (\ref{Energy-identity}). \qed\\

Now we want to estimate $\|\nabla u(t)\|_{H^1}^2$ in terms of $\Phi(t)$.
\begin{lemma}\label{le:2.2} There exists  $M\in\mathcal F$ such that for $0<t<T$, it holds
 \beq\label{H^2 of u}
\begin{split}
\|\nabla u(t)\|_{H^1}\le &M(\Phi(t)).
\end{split}
\eeq
\end{lemma} 
\pf By the standard $H^2$-estimate of the Lam\'e equation with respect to the boundary condition (\ref{clcboundary1}) or (\ref{clcboundary2})
or (\ref{clcboundary3}), (\ref{infty_bound_rho}), and H\"older's
inequality, we have \beq\label{H^2 of u:1}
\begin{split}
\|\nabla u\|_{H^1}^2\lesssim& \|\mathcal L u\|_{L^2}^2+\|\nabla u\|_{L^2}^2\\
\les &\|\rho u_t\|_{L^2}^2+\|\rho u\cdot\nabla u\|_{L^2}^2+\|\nabla (P(\rho))\|_{L^2}^2+\|\de d\cdot\nabla d\|_{L^2}^2+\|\nabla u\|_{L^2}^2\\
\les &\|\rho\|_{L^{\infty}}
\|\sqrt{\rho}u_t\|_{L^2}^2+\|\rho\|_{L^{\infty}}^2\|
u\|_{L^6}^2\|\nabla u\|_{L^3}^2
+B_P^2(\|\rho\|_{L^{\infty}})\|\nabla \rho\|_{L^2}^2\\&+\|\de
d\|_{L^3}^2\|\nabla d\|_{L^6}^2+\|\nabla u\|_{L^2}^2\\
\le& M(\Phi)(1+\|u\|_{L^6}^2\|\nabla u\|_{L^3}^2)+C\|\de
d\|_{L^3}^2\|\nabla d\|_{L^6}^2
\end{split}
\eeq 
for some $M\in\mathcal F$.
By the interpolation inequality, Sobolev's inequality\footnote{ when $\Omega=B_R$ for $R\ge 1$, by simple scalings, one has
$$\|f\|_{L^6(B_R)}\le C\left(R^{-1}\|f\|_{L^2(B_R)}+\|\nabla f\|_{L^2(B_R)}\right)
\le C\|f\|_{H^1(B_R)}.$$},  we obtain
\beq\label{H^2 of u:1.3}
\begin{split}
\| u\|_{L^6}^2\|\nabla u\|_{L^3}^2
\leq& C\|\nabla u\|_{L^2}^3\|\nabla u\|_{H^1}.
\end{split}
\eeq 
Similar to (\ref{H^2 of u:1.3}), by
(\ref{Energy-identity}), we obtain
\beq\label{H^2 of u:1.4}
\begin{split}
&\|\de d\|_{L^3}^2\|\nabla d\|_{L^6}^2
\les\|\de d\|_{L^2}\|\de d\|_{L^6}\|\nabla d\|_{H^1}^2\\
\les& \|\de d\|_{H^1}^2\|\nabla d\|_{L^2}^2+\|\de d\|_{H^1}^2\|\nabla^2 d\|_{L^2}^2
\les M(\Phi)
\end{split}
\eeq for some $M\in\mathcal F$.
Substituting (\ref{H^2 of u:1.3}), (\ref{H^2 of u:1.4})  into (\ref{H^2 of u:1}),  and using  (\ref{Energy-identity}) and Cauchy's inequality, we have 
\bex
\begin{split}
\|\nabla u\|_{H^1}^2\le&\frac{1}{2}\|\nabla
u\|_{H^1}^2+M(\Phi(t))
\end{split}
\eex for some $M\in\mathcal F$.
This gives (\ref{H^2 of u}) and completes the proof. 
\endpf\\

Now we want to estimate $\|\sqrt{\rho} u_t\|_{L^2}$. More precisely, we have
\begin{lemma}\label{le:nabla u_t} There exists $M\in\mathcal F$ such tha for any $0<t<T$, it holds
\beq\label{rho u_tt estimates 3} \int_\Omega\rho
|u_t|^2dx+\int_0^t\int_\Omega|\nabla u_t|^2dxds \leq C \mathcal M(\rho_0,
u_0, d_0)+ \int_0^t M(\Phi(s))\,ds.\eeq
\end{lemma}

\pf Differentiating (\ref{clc-2}) with respect to $t$, we have\footnote{here we have used the fact that
$\Delta d\cdot\nabla d=\nabla\cdot(\nabla d\otimes\nabla d-\frac{1}{2}|\nabla d|^2 \mathbb I_3),$
where $\nabla
d\otimes\nabla d=\left(d_{x_i}\cdot d_{x_j}\right)_{1\le i, j\le 3}$ and  $\mathbb I_3$ is the identity matrix of 
order $3$.}
\beq\label{rho u_tt}\begin{split} & \rho u_{tt}+\rho
u\cdot\nabla u_t+\rho_tu_t+\rho_tu\cdot\nabla u+\rho u_t\cdot\nabla u+\nabla (P(\rho))_t\\
=&(2\mu+\lambda)\nabla\di u_t-\mu\nabla\times(\nabla\times
u_t)-\nabla\cdot(\nabla d_t\otimes\nabla d+\nabla d\otimes\nabla
d_t-\nabla d\cdot\nabla d_t\ \mathbb I_3).
\end{split}
\eeq 
Multiplying (\ref{rho u_tt}) by $u_t$,
integrating the resulting equations over $\Omega$, and using
(\ref{clc-1}) and integration by parts, we have \beq\label{2.9}
\begin{split}
&\frac{1}{2}\frac{d}{dt}\int\rho
|u_t|^2\,dx+\int\left((2\mu+\lambda)|\di
u_t|^2+\mu|\nabla\times u_t|^2\right)\,dx\\
=&-2\int\rho u u_t\cdot\nabla
u_t\,dx-\int\rho_tu\cdot\nabla u\cdot u_t\,dx-\int\rho
u_t\cdot\nabla u\cdot
u_t\,dx+\int P'(\rho)\rho_t\di
u_t\,dx\\&+\int (\nabla d_t\otimes\nabla d+\nabla d\otimes\nabla
d_t-\nabla d\cdot\nabla d_t\ \mathbb I_3):\nabla
u_t\,dx=\sum\limits_{i=1}^5II_i.
\end{split}
\eeq 
By H\"older's inequality, Sobolev's inequality,   (\ref{infty_bound_rho}),  and (\ref{H^2 of u}),
we have \beq\label{2.9:1}
\begin{split}
|II_1|\les&\|\nabla
u_t\|_{L^2}\|\sqrt{\rho}u_t\|_{L^2}\|\sqrt{\rho}u\|_{L^\infty}\\
\les&\|\nabla
u_t\|_{L^2}\|\sqrt{\rho}u_t\|_{L^2}\|\sqrt{\rho}\|_{L^\infty}\|\nabla
u\|_{H^1}
\les  M(\Phi)\|\nabla u_t\|_{L^2}
\end{split}
\eeq
for some $M\in\mathcal F$.

By (\ref{clc-1}), H\"older's inequality, Sobolev's inequality,
(\ref{infty_bound_rho}), and (\ref{H^2 of u}), we have 
\beq\label{2.9:2}
\begin{split}
|II_2| =&\left|\int\rho u\cdot\nabla(u\cdot\nabla u\cdot
u_t)\,dx\right|\\=&
\left|\int\rho u\cdot\left(\nabla u\cdot\nabla u\cdot u_t+u\cdot\nabla\nabla u\cdot u_t+u\cdot\nabla u\cdot \nabla u_t\right)\,dx\right|\\
\les& \|\rho\|_{L^\infty}\|u\|_{L^6}^2\|\nabla u\|_{L^6}\|\nabla
u_t\|_{L^2}+\|\rho\|_{L^\infty}\|u\|_{L^6}^2\|\nabla^2u\|_{L^2}\|u_t\|_{L^6}\\
&+\|\sqrt{\rho}\|_{L^\infty}\|u\|_{L^6}\|\nabla
u\|_{L^6}^2\|\sqrt{\rho}u_t\|_{L^2}\\
\les& M(\Phi)(1+\|\nabla u_t\|_{L^2})
\end{split}
\eeq for some $M\in\mathcal F$.
For $II_3$, by (\ref{H^2 of u}) we have \beq\label{2.9:3}
\begin{split}
|II_3|\les&
\|\sqrt{\rho}\|_{L^\infty}\|\sqrt{\rho}u_t\|_{L^2}\|\nabla
u\|_{L^3}\|u_t\|_{L^6}\\
\les& \|\sqrt{\rho}\|_{L^\infty}\|\sqrt{\rho}u_t\|_{L^2}   \|\nabla
u\|_{L^2}^{\frac{1}{2}}\|\nabla
u\|_{H^1}^{\frac{1}{2}}\|\nabla u_t\|_{L^2}
\les M(\Phi)\|\nabla u_t\|_{L^2}
\end{split}
\eeq 
for some $M\in\mathcal F$.
For $II_4$, by (\ref{clc-1}), (\ref{infty_bound_rho}), and (\ref{H^2 of u}) we have \beq\label{2.9:4}
\begin{split}
|II_4|\les& B_P(\|\rho\|_{L^\infty})\|\rho_t\|_{L^2}\|\di
u_t\|_{L^2}\\
\les&B_P(\|\rho\|_{L^\infty}) (\|\nabla\rho\|_{L^2}\|u\|_{L^\infty}+\|\rho\|_{L^\infty}\|\di
u\|_{L^2})\|\di u_t\|_{L^2}\\
\les & M(\Phi)\|\nabla u_t\|_{L^2}
\end{split}
\eeq 
for some $M\in\mathcal F$.
For $II_5$, by (\ref{Energy-identity}) we have \beq\label{2.9:5}
\begin{split}
|II_5|\les& \int_\Omega |\nabla d||\nabla d_t||\nabla u_t|dx 
\les\|\nabla u_t\|_{L^2}\|\nabla d\|_{L^\infty}\|\nabla
d_t\|_{L^2}\\
\les&\|\nabla u_t\|_{L^2}\|\nabla d\|_{H^2}\|\nabla d_t\|_{L^2}\\
\les&\|\nabla u_t\|_{L^2}(\|\nabla d\|_{L^2}+\|\nabla^2 d\|_{H^1})\|\nabla d_t\|_{L^2}
\leq(C+ M(\Phi))\|\nabla u_t\|_{L^2}\|\nabla d_t\|_{L^2}
\end{split}
\eeq
for some $M\in\mathcal F$.
Substituting (\ref{2.9:1})-(\ref{2.9:5}) into (\ref{2.9}), and using Cauchy's inequality, we have
\beq\label{2.15}
\begin{split}
&\frac{1}{2}\frac{d}{dt}\int\rho
|u_t|^2dx+\frac{1}{C}\int|\nabla
u_t|^2dx \\
\le&\frac{1}{2C}\int|\nabla u_t|^2dx+M(\Phi)+(C+M(\Phi))\|\nabla d_t\|_{L^2}^2
\end{split}
 \eeq 
for some $M\in\mathcal F$, 
where we have used the following inequality due to \cite{Wolf}:
if (i) either $\Omega$ is simply connected and $u\cdot\nu=0$ on $\partial\Omega$ or (ii) $u=0$ on $\partial\Omega$\ \footnote{in fact, in this case, the inequality (\ref{div-curl-ineq}) is an equality.},
then
\beq \label{div-curl-ineq} \|\nabla u_t\|_{L^2}\les\|\mathrm{div}
u_t\|_{L^2}+\|\nabla\times u_t\|_{L^2}. \eeq 
By (\ref{2.15}), we
have \beq\label{rho u_tt estimates 1}
\begin{split}
\frac{d}{dt}\int\rho |u_t|^2dx+\frac{1}{C}\int|\nabla
u_t|^2dx\les M(\Phi)
+(C+M(\Phi))\|\nabla d_t\|_{L^2}^2.
\end{split}
\eeq Differentiating (\ref{clc-3}) with respect to $x$, we have
\beq\label{exss8.1} \nabla d_t-\nabla\de d=\nabla(|\nabla
d|^2d)-\nabla(u\cdot\nabla d). \eeq
 From (\ref{exss8.1}), we have\ \footnote{here we also use the Sobolev's inequality:
$\|\nabla d\|_{L^\infty(\Omega)}\le C\|\nabla d\|_{H^2(\Omega)}$
and the fact that $C$ can be chosen independent of $R$ when $\Omega=B_R$ for $R\ge 1$.}
 \beq\label{nabla d_t}
\begin{split}
\|\nabla d_t\|_{L^2}\lesssim&\|\nabla u\cdot\nabla
d\|_{L^2}+\|u\cdot\nabla^2d\|_{L^2}+\|\nabla\Delta
d\|_{L^2}+\|\nabla d\|_{L^6}^3+\|\nabla d\cdot\nabla^2d\|_{L^2}\\
\lesssim&\|\nabla d\|_{L^\infty}\|\nabla
u\|_{L^2}+\|u\|_{L^6}\|\nabla^2d\|_{L^3}+\|\nabla\Delta
d\|_{L^2}+(1+\|\nabla^2 d\|_{L^2})^3\\&+\|\nabla
d\|_{L^\infty}\|\nabla^2d\|_{L^2}\\ \lesssim&\|\nabla
d\|_{H^2}\|\nabla u\|_{L^2}+\|\nabla
u\|_{L^2}\|\nabla^2d\|_{H^1}+\|\nabla\Delta d\|_{L^2}+(1+\|\nabla^2
d\|_{L^2})^3\\&+\|\nabla d\|_{H^2}\|\nabla^2d\|_{L^2}\\
\lesssim& M(\Phi)+1
\end{split}
\eeq for some $M\in \mathcal F$.

 Substituting (\ref{nabla d_t}) into
(\ref{rho u_tt estimates 1}), and using Cauchy's inequality, we have
\beq\label{rho u_tt estimates 2}
\begin{split}
\frac{d}{dt}\int\rho |u_t|^2dx+\frac{1}{C}\int|\nabla
u_t|^2dx \leq M(\Phi)+C
\end{split}
\eeq for some $M\in\mathcal F$.
Integrating (\ref{rho u_tt estimates 2}) over $(0,t)$, and
using (\ref{clc-2}), and (\ref{first3.1}), we have \bex
\int\rho |u_t|^2dx+\int_0^t\int_\Omega|\nabla u_t|^2dxds
\leq C\int\rho|u_t|^2\,dx\Big|_{t=0}+\int_0^t M(\Phi(s))ds+C\\
\leq C\mathcal M(\rho_0, u_0, d_0)+ \int_0^t M(\Phi(s))ds
\eex 
for some $M\in\mathcal F$. 
This completes the proof.
\endpf\\

As an immediate consequence of Lemma \ref{le:nabla u_t},  we obtain an estimate of $\|\nabla u\|_{L^2}$.
\begin{lemma}\label{le:2.4} There exists $M\in\mathcal F$ such that for $0<t<T$, it holds
\beq\label{H^1 of u:2} 
\begin{split}
\int |\nabla u(t)|^2\,dx \leq C\mathcal M(\rho_0, u_0, d_0)+
\int_0^t M(\Phi(s))\,ds.
\end{split}
\eeq
\end{lemma}

\pf By Cauchy's inequality, Lemma \ref{le:2.1}), Lemma \ref{le:2.2}, and
Lemma \ref{le:nabla u_t}, we have
\bex\begin{split} \int|\nabla u|^2(t)dx=&\int|\nabla
u_0|^2dx+2\int_0^t\int_\Omega\nabla u\cdot\nabla u_tdxds\\
\leq&C+\int_0^t\int_\Omega|\nabla
u|^2dxds+\int_0^t\int_\Omega|\nabla u_t|^2dxds\\ \leq& C\mathcal M(\rho_0, u_0,
d_0)+ \int_0^tM(\Phi(s))\,ds
\end{split}
\eex
for some $M\in\mathcal F$. 
This completes the proof. \endpf

\begin{lemma}\label{nabla rho est} There exists $M\in\mathcal F$ such that for $0<t<T$, it holds
\beq\label{w^1,q of rho:2}
\begin{split}
\|\rho(t)\|_{H^1\cap W^{1,q}}\leq \exp\left\{C\mathcal M(\rho_0, u_0,
d_0)+ C\int_0^t M(\Phi(s))\,ds\right\}.
\end{split}
\eeq
\end{lemma}
\pf It follows from \cite{CCK} (page 249,  (2.11)) that \beq\label{w^1,q of
rho:1}
\begin{split}
\|\rho(t)\|_{H^1\cap W^{1,q}}\le \|\rho_0\|_{H^1\cap
W^{1,q}}\exp\left\{C\int_0^t\|\nabla u\|_{H^1\cap
D^{1,q}}ds\right\}.
\end{split}
\eeq By $W^{2,q}$-estimate of the Lam\'e equation under either Dirichlet boundary condition (\ref{clcboundary2})
or the Navier-slip boundary condition (\ref{clcboundary3}) (see Lemma {\ref{2p-est-lemma}} below),
(\ref{clc-2}), and Sobolev's inequality, we have
 \beq\label{W^2,q of u:1}
\begin{split}
\|\nabla^2u\|_{L^q}\lesssim&\|\rho u_t\|_{L^q}+\|\rho u\cdot\nabla
u\|_{L^q}+\|\nabla (P(\rho))\|_{L^q}+\|\nabla d\cdot \Delta
d\|_{L^q}=\sum\limits_{i=1}^4III_i.
\end{split}
\eeq If $q=6$, then by Sobolev's inequality we have \beq\label{W^2,q of u:1.1.1}
\begin{split}
III_1\lesssim\|\rho\|_{L^\infty}\|u_t\|_{L^6}\les \Phi\|\nabla
u_t\|_{L^2}.
\end{split}
\eeq 
If $q\in(3,6)$, then by H\"older's inequality and Sobolev's inequality, we have 
\beq\label{W^2,q of u:1.1.2}
\begin{split}
III_1\lesssim\|\rho\|_{L^\frac{6q}{6-q}}\|u_t\|_{L^6}
\les\|\rho\|_{L^1}^{\frac{6-q}{6q}}\|\rho\|_{L^\infty}^{1-\frac{6-q}{6q}}\|\nabla u_t\|_{L^2}
\les \Phi\|\nabla u_t\|_{L^2},
\end{split}
\eeq where we have used the fact that $\int\rho
dx=\int\rho_0dx$. From (\ref{W^2,q of u:1.1.1}) and (\ref{W^2,q of u:1.1.2}), we have
that for $q\in(3,6]$,
 \beq\label{W^2,q of u:1.1}
\begin{split}
III_1\lesssim \Phi\|\nabla u_t\|_{L^2}.
\end{split}
\eeq For $III_2$,  if $q\in (3,6]$, then by similar arguments, Lemma 2.2, and Lemma 2.3,  we have
\beq\label{W^2,q of u:1.2}
\begin{split}
III_2\lesssim\Phi\|\nabla u\|_{H^1}^2
\leq M(\Phi)
\end{split}
\eeq
for some $M\in\mathcal F$. For $III_3$ and $III_4$, if $q\in(3,6]$, then we have
\beq\label{W^2,q of u:1.3, 1.4}
\begin{split}
III_3+III_4\leq CB_P(\|\rho\|_{L^\infty})\|\nabla\rho\|_{L^q}+\|\nabla d\|_{H^2}^2
\leq M(\Phi)
\end{split}
\eeq
for some $M\in\mathcal F$.
 Substituting (\ref{W^2,q of u:1.1}), (\ref{W^2,q of u:1.2}) and
(\ref{W^2,q of u:1.3, 1.4}) into (\ref{W^2,q of u:1}), we have
\beq\label{W^2,q of u}
\begin{split}
\|\nabla^2u\|_{L^q}\lesssim&\Phi\|\nabla
u_t\|_{L^2}+M(\Phi)\le \|\nabla u_t\|_{L^2}^2+M(\Phi)
\end{split}
\eeq 
for some $M\in\mathcal F$.
Integrating (\ref{W^2,q of u}) over $(0,t)$, and using Cauchy's
inequality and (\ref{rho u_tt estimates 3}), we have
\beq\label{W^2,q of u:2}
\begin{split}
\int_0^t\|\nabla^2u\|_{L^q}\leq C\mathcal M(\rho_0, u_0, d_0)+
\int_0^tM(\Phi(s))\,ds.
\end{split}
\eeq
 Substituting (\ref{H^2 of u}) and (\ref{W^2,q of u:2}) into (\ref{w^1,q of rho:1}), we have
 \bex
\begin{split}
\|\rho(t)\|_{H^1\cap W^{1,q}}\lesssim \exp\left\{C\mathcal M(\rho_0, u_0, d_0)+
C\int_0^t M(\Phi(s))\,ds\right\}
\end{split}
\eex for some $M\in\mathcal F$.
This completes the proof.
\endpf

\begin{lemma}\label{le:H^2 of d} There exists $M\in\mathcal F$ such that for any $0<t<T$, it holds
\beq\label{exss8.6}
\begin{split}
\|\nabla^2 d\|^{2}_{L^2}+\int_0^t\|\nabla d_t\|^2_{L^2}ds\leq C+
\int_0^t M(\Phi(s))\,ds.
\end{split}
\eeq
\end{lemma}
\pf Multiplying (\ref{exss8.1}) by $\nabla d_t$ and integrating
over $\om$, using integration by parts and $\frac{\partial
d_t}{\partial\nu}=0$ on $\partial\om$, we obtain \bex
\begin{split}
\|\nabla d_t\|^2_{L^2}+\frac{1}{2}\frac{d}{dt}\|\de d\|^{2}_{L^2}=&\int\left[\nabla(|\nabla d|^2d)-\nabla(u\cdot\nabla d)\right]\nabla d_t\,dx\\
\leq&\frac{1}{2}\|\nabla d_t\|^2_{L^2} +C\int|\nabla(|\nabla
d|^2d)|^2\,dx+C\int|\nabla(u\cdot\nabla d)|^2\,dx.
\end{split}
\eex 
Thus we have
\beq\label{exss8.2} \|\nabla d_t\|^2_{L^2}+\frac{d}{dt}\|\de
d\|^{2}_{L^2}\les \int|\nabla(|\nabla
d|^2d)|^2\,dx+\int|\nabla(u\cdot\nabla d)|^2\,dx. \eeq 
Similar to the proof of (\ref{nabla d_t}), we obtain 
\beq\label{exss8.5} \|\nabla
d_t\|^2_{L^2}+\frac{d}{dt}\|\de d\|^{2}_{L^2}\leq M(\Phi) \eeq
for some $M\in\mathcal F$.
Integrating (\ref{exss8.5}) over $(0,t)$ and applying $W^{2,2}$-estimate of the equation (\ref{clc-3}), we have \bex
\begin{split}
\|\nabla^2 d\|^{2}_{L^2}+\int_0^t\|\nabla
d_t\|^2_{L^2}ds\leq\|\nabla^2 d_0\|^{2}_{L^2}+\int_0^t M(\Phi(s))\,ds\leq
&C+\int_0^t M(\Phi(s))\,ds.
\end{split}
\eex 
This completes the proof.
\endpf

\begin{lemma}\label{le:H^3 of d} There exists $M\in\mathcal F$ such that for $0<t<T$, it holds
\beq\label{H^3 of d}
\begin{split}
\|\nabla^3d\|^{2}_{L^2}+\int_0^t\|\nabla^2 d_t\|^2_{L^2}\,ds
\leq \left(C\mathcal M(\rho_0, u_0, d_0)+\int_0^t M(\Phi(s))\,ds\right)^4.
\end{split}
\eeq
\end{lemma}
\pf Multiplying (\ref{exss8.1}) by $\nabla\de d_t$, integrating
over $\om$, using $\frac{\partial
d_t}{\partial\nu}=0$ on $\partial\om$  and integration by parts,
we obtain \beq\label{exss9.1}
\begin{split}
\|\de d_t\|^2_{L^2}+\frac{1}{2}\frac{d}{dt}\|\nabla\de
d\|^{2}_{L^2}=&\int\left[\nabla(u\cdot\nabla d)-\nabla(|\nabla
d|^2d)\right]\cdot\nabla\de d_t\,dx
\\=&\frac{d}{dt}\int \left[\nabla(u\cdot\nabla d)-\nabla(|\nabla d|^2d)\right]\cdot\nabla\de d\,dx\\
-&\int\frac{\partial}{\partial t}\left[\nabla(u\cdot\nabla d)-\nabla(|\nabla
d|^2d)\right]\cdot\nabla\de d\,dx.
\end{split}
\eeq  Now we need to estimate the
second term of right side as follows.
\beq\label{H^3 of d: 1}
\begin{split}
-\int \frac{\partial}{\partial t}[\nabla(u\cdot\nabla d)]\cdot\nabla\de d\,dx
=&-\int[\nabla u_t\cdot\nabla d+\nabla u\cdot\nabla d_t+
u_t\cdot\nabla^2 d+u\cdot\nabla^2 d_t]\cdot\nabla\de d\,dx\\
=&\sum\ls_{i=1}^4 IV_{i}.
\end{split}
\eeq By H\"older's inequality and Sobolev's inequality, we have
\beq\label{H^3 of d: 1.1}
\begin{split}
|IV_1|\les\|\nabla u_t\|_{L^2}\|\nabla d\|_{L^\infty}\|\nabla\de d\|_{L^2}
\les\|\nabla u_t\|_{L^2}\|\nabla d\|_{H^2}^2
\les M(\Phi)+\|\nabla u_t\|_{L^2}^2
\end{split}
\eeq 
for some $M\in\mathcal F$.

By H\"older's inequality, Sobolev's inequality, (\ref{H^2 of u}), (\ref{nabla d_t}) and Young's inequality, we obtain 
\beq\label{H^3 of d: 1.2}
\begin{split}
|IV_2|\les&\|\nabla u\|_{L^6}\|\nabla d_t\|_{L^3}\|\nabla\de d\|_{L^2}\\
\les&\|\nabla u\|_{H^1}\|\nabla d_t\|_{H^1}\|\nabla\de d\|_{L^2}\\
\les&\|\nabla u\|_{H^1}\|\nabla^2 d_t\|_{L^2}\|\nabla\de d\|_{L^2}+\|\nabla u\|_{H^1}\|\nabla d_t\|_{L^2}\|\nabla\de d\|_{L^2}\\
\le& \va\|\nabla^2 d_t\|^2_{L^2}+M(\Phi)
\end{split}
\eeq 
for some $M\in\mathcal F$.

By H\"older's inequality, Sobolev's inequality and Cauchy's inequality, we obtain
 \beq\label{H^3 of d: 1.3}
\begin{split}
|IV_3|\les&\|u_t\|_{L^6}\|\nabla^2 d\|_{L^3}\|\nabla\de d\|_{L^2}\\
\les&\|\nabla u_t\|_{L^2}\|\nabla^2 d\|_{H^1}\|\nabla\de d\|_{L^2}
\les M(\Phi)+\|\nabla u_t\|_{L^2}^2
\end{split}
\eeq 
for some $M\in\mathcal F$.

By H\"older's inequality, Sobolev's inequality, (\ref{H^2 of u}) and Cauchy's inequality, we obtain
 \beq\label{H^3 of d: 1.4}
\begin{split}
|IV_4|\les&\|u\|_{L^\infty}\|\nabla^2 d_t\|_{L^2}\|\nabla\de d\|_{L^2}
\les\|\nabla u\|_{H^1}\|\nabla^2 d_t\|_{L^2}\|\nabla\de d\|_{L^2}\\
\le& \va\|\nabla^2 d_t\|^2_{L^2}+M(\Phi)
\end{split}
\eeq 
for some $M\in\mathcal F$.

 Combining (\ref{H^3 of d: 1}), (\ref{H^3 of d: 1.1}), (\ref{H^3 of d:
1.2}), (\ref{H^3 of d: 1.3}) and (\ref{H^3 of d: 1.4}), we obtain
\beq\label{exss9.2}
\begin{split}
-\int \frac{\partial}{\partial t}[\nabla(u\cdot\nabla d)]\cdot\nabla\de d\,dx
\le&2\va\|\nabla^2 d_t\|^2_{L^2}+C\|\nabla
u_t\|_{L^2}^2+ M(\Phi)
\end{split}
\eeq 
for some $M\in\mathcal F$.

By Leibniz's rule and the fact $|d|=1$, we
have \beq\label{H^3 of d: 2}
\begin{split}
&\int\frac{\partial}{\partial t}[\nabla(|\nabla d|^2d)]\cdot\nabla\de d\,dx\\
\les&\int [|\nabla d|^2|\nabla d_t|+|\nabla d_t||\nabla^2
d|+|\nabla d||\nabla^2 d_t|+|\nabla d||\nabla^2 d||d_t|]|\nabla\de
d|\,dx
=\sum\ls_{i=1}^4 V_{i}.
\end{split}
\eeq By H$\ddot{\mbox{o}}$lder's inequality, Sobolev's
inequality and (\ref{nabla d_t}),  Cauchy inequality, and Young inequality, we obtain 
\beq\label{H^3 of d: 2.1}
\begin{split}
|V_1|\les&\|\nabla d\|_{L^\infty}^2\|\nabla d_t\|_{L^2}\|\nabla\de d\|_{L^2}
\les\|\nabla d\|_{H^2}^2\|\nabla d_t\|_{L^2}\|\nabla\de d\|_{L^2}
\leq M(\Phi),
\end{split}
\eeq 
\beq\label{H^3 of d: 2.2}
\begin{split}
|V_2|\les&\|\nabla d_t\|_{L^6}\|\nabla^2 d\|_{L^3}\|\nabla\de d\|_{L^2}
\les\|\nabla d_t\|_{H^1}\|\nabla^2 d\|_{H^1}^2\\
\les&\Phi(\|\nabla^2 d_t\|_{L^2}+\|\nabla d_t\|_{L^2})
\le \va \|\nabla^2 d_t\|^2_{L^2}+M(\Phi),
\end{split}
\eeq 
\beq\label{H^3 of d: 2.3}
\begin{split}
|V_3|\les&\|\nabla d\|_{L^\infty}\|\nabla^2 d_t\|_{L^2}\|\nabla\de d\|_{L^2}\\
\les&\|\nabla d\|_{H^2}\|\nabla^2 d_t\|_{L^2}\|\nabla\de d\|_{L^2}
\le\varepsilon\|\nabla^2 d_t\|_{L^2}^2+M(\Phi),
\end{split}
\eeq
 \beq\label{H^3 of d: 2.4}
\begin{split}
|V_4|\les&\|d_t\|_{L^6}\|\nabla d\|_{L^\infty}\|\nabla^2 d\|_{L^3}\|\nabla\de d\|_{L^2}\\
\les&\|d_t\|_{H^1}\|\nabla d\|_{H^2}\|\nabla^2 d\|_{H^1}\|\nabla\de d\|_{L^2}
\les\|d_t\|_{H^1} M(\Phi)
\end{split}
\eeq 
for some $M\in\mathcal F$. Notice that
\beq\label{L^2 of d_t}
\begin{split}
\|d_t\|_{L^2}\lesssim&\|\Delta d\|_{L^2}+\|\nabla
d\|_{L^4}^2+\|u\cdot\nabla d\|_{L^2}\\ \lesssim& \|\nabla
d\|_{H^1}^2+\|u\|_{L^6}\|\nabla d\|_{L^3}+1 \lesssim\|\nabla
d\|_{H^1}^2+\|\nabla u\|_{L^2}\|\nabla d\|_{H^1}+1
\lesssim\Phi.
\end{split}
\eeq 
Thus by (\ref{nabla d_t}), (\ref{H^3 of d: 2.4}) and (\ref{L^2 of d_t}),
we have 
\beq\label{H^3 of d: 2.5}
\begin{split}
|V_4|\leq M(\Phi)
\end{split}
\eeq 
for some $M\in\mathcal F$.
Combining (\ref{H^3 of d: 2}), (\ref{H^3 of d: 2.1}), (\ref{H^3 of d:
2.2}), (\ref{H^3 of d: 2.3}) and (\ref{H^3 of d: 2.5}), we have
\beq\label{H^3 of d: 2 total}
\begin{split}
\int\frac{\partial}{\partial t}[\nabla(|\nabla d|^2d)]\cdot\nabla\de d\,dx \le
2\varepsilon\|\nabla^2 d_t\|_{L^2}^2+M(\Phi)
\end{split}
\eeq
for some $M\in\mathcal F$.

 Putting (\ref{exss9.2}) and (\ref{H^3 of d: 2 total}) into
(\ref{exss9.1}), we obtain \beq\label{exss9.4}
\begin{split}
\|\de d_t\|^2_{L^2}+\frac{1}{2}\frac{d}{dt}\|\nabla\de
d\|^{2}_{L^2}\le&\frac{d}{dt}\int\left[\nabla(u\cdot\nabla
d)-\nabla(|\nabla d|^2d)\right]
\cdot\nabla\de d\,dx\\
+&4\va\|\nabla^2 d_t\|^2_{L^2}+C\|\nabla
u_t\|_{L^2}^2+M(\Phi)
\end{split}
\eeq 
for some $M\in\mathcal F$.
Integrating (\ref{exss9.4}) over $(0,t)$, using $H^k$ ($k=2,3$)
estimate of the elliptic equations, and choosing $\varepsilon$
small enough, we have \beq\label{exss9.5}
\begin{split}
&\|\nabla^3d\|^{2}_{L^2}+\int_0^t\|\nabla^2 d_t\|^2_{L^2}ds\\
\les&\int\left|\nabla(u\cdot\nabla d)-\nabla(|\nabla
d|^2d)\right| |\nabla\de d|\,dx+\int\left|\nabla(u_0\cdot\nabla d_0)-\nabla(|\nabla
d_0|^2d_0)\right| |\nabla\de d_0|\,dx\\
&+\|\nabla^3d_0\|^{2}_{L^2}+\int_0^t\|\nabla
u_t\|_{L^2}^2ds+\int_0^t M(\Phi(s))\,ds.
\end{split}
\eeq
For the first term of right side of (\ref{exss9.5}), we have
\beq\label{exss9.8}
\begin{split}
&\int\left|\nabla(u\cdot\nabla d)-\nabla(|\nabla
d|^2d)\right| |\nabla\de d|\,dx\\
\les&\int\left(|\nabla u||\nabla d|+| u||\nabla^2 d|+|\nabla d|^3+|\nabla d||\nabla^2 d|\right) |\nabla\de d|\,dx
=\sum\ls_{i=1}^4VI_i.
\end{split}
\eeq
By H$\ddot{\mbox{o}}$lder's inequality, Nirenberg's interpolation inequality, (\ref{Energy-identity}),
and Young's inequality, 
we obtain
\beq\label{exss9.9}
\begin{split}
|VI_1|\les&\|\nabla d\|_{L^\infty}\|\nabla u\|_{L^2}\|\nabla\Delta
d\|_{L^2}
\les\|\nabla d\|_{L^2}^{\frac{1}{4}}\|\nabla
d\|_{H^2}^{\frac{3}{4}}\|\nabla u\|_{L^2}\|\nabla\Delta
d\|_{L^2}\\
\les&\|\nabla d\|_{H^1}^{\frac{3}{4}}\|\nabla
u\|_{L^2}\|\nabla^3
d\|_{L^2}+\|\nabla^3d\|_{L^2}^{\frac{7}{4}}\|\nabla
u\|_{L^2}\\
\leq&\va\|\nabla^3d\|_{L^2}^2+C(\|\nabla d\|_{H^1}^{\frac{3}{2}}\|\nabla
u\|_{L^2}^2+\|\nabla u\|_{L^2}^8),
\end{split}
\eeq 
\beq\label{exss9.10}
\begin{split}
|VI_2|\les&\|u\|_{L^6}\|\nabla^2 d\|_{L^3}\|\nabla\Delta
d\|_{L^2}
\les\|\nabla u\|_{L^2}\|\nabla^2
d\|_{L^2}^{\frac{1}{2}}\|\nabla^2
d\|_{H^1}^{\frac{1}{2}}\|\nabla^3d\|_{L^2}\\
\les&\|\nabla u\|_{L^2}\|\nabla^2 d\|_{L^2}\|\nabla^3
d\|_{L^2}+\|\nabla u\|_{L^2}\|\nabla^2
d\|_{L^2}^{\frac{1}{2}}\|\nabla^3
d\|_{L^2}^{\frac{3}{2}}\\
\leq&\va\|\nabla^3d\|_{L^2}^2+C\|\nabla^2 d\|_{L^2}^2(\|\nabla u\|_{L^2}^{2}+\|\nabla u\|_{L^2}^4),
\end{split}
\eeq
\beq\label{exss9.11}
\begin{split}
|VI_3|\les&\|\nabla d\|_{L^6}^3\|\nabla\Delta d\|_{L^2}
\les\|\nabla d\|_{H^1}^3\|\nabla^3 d\|_{L^2}
\leq\va\|\nabla^3d\|_{L^2}^2+C\|\nabla d\|_{H^1}^6,
\end{split}
\eeq
and
\beq\label{exss9.12}
\begin{split}
|VI_4|\les&\|\nabla
d\|_{L^\infty}\|\nabla^2 d\|_{L^2}\|\nabla\Delta d\|_{L^2}
\les\|\nabla
d\|_{L^2}^{\frac{1}{4}}\|\nabla d\|_{H^2}^{\frac{3}{4}}\|\nabla^2
d\|_{L^2}\|\nabla^3 d\|_{L^2}\\
\les&\|\nabla d\|_{H^1}^{\frac{3}{4}}\|\nabla^2d\|_{L^2}\|\nabla^3
d\|_{L^2}+\|\nabla^3d\|_{L^2}^{\frac{7}{4}}\|\nabla^2d\|_{L^2}\\
\leq&\va\|\nabla^3d\|_{L^2}^2+C(\|\nabla d\|_{H^1}^{\frac{3}{2}}\|\nabla^2d\|_{L^2}^2+\|\nabla^2d\|_{L^2}^8).
\end{split}
\eeq
Combining (\ref{exss9.8}), (\ref{exss9.9}), (\ref{exss9.10}), (\ref{exss9.11}) and (\ref{exss9.12}), we obtain
\beq\label{exss9.13}
\begin{split}
&\int\left|\nabla(u\cdot\nabla d)-\nabla(|\nabla
d|^2d)\right| |\nabla\de d|\,dx\\
\leq&4\va\|\nabla^3d\|_{L^2}^2+C\|\nabla d\|_{H^1}^{\frac{3}{2}}(\|\nabla u\|_{L^2}^2+\|\nabla^2d\|_{L^2}^2)
+C\|\nabla^2 d\|_{L^2}^2(\|\nabla
u\|_{L^2}^2+\|\nabla u\|_{L^2}^4)\\
+&C(\|\nabla d\|_{H^1}^6+\|\nabla^2d\|_{L^2}^8+\|\nabla u\|_{L^2}^8)\\
\leq &4\va\|\nabla^3d\|_{L^2}^2+\left(C\mathcal M(\rho_0, u_0, d_0)+\int_0^t M(\Phi(s))\,ds\right)^4
\end{split}
\eeq 
for some $M\in\mathcal F$, where we have used Lemma 2.2, Lemma 2.5, and Lemma 2.7 in the last step.

Substituting (\ref{exss9.13}) into (\ref{exss9.5}), choosing $\va$ small enough,
and using (\ref{rho u_tt estimates 3}), Cauchy's inequality, Lemma \ref{le:2.4} and
(\ref{exss8.6}), we have
\bex
\begin{split}
&\|\nabla^3d\|^{2}_{L^2}+\int_0^t\|\nabla^2 d_t\|^2_{L^2}ds
\leq \left(C\mathcal M(\rho_0, u_0, d_0)+\int_0^t M(\Phi(s))\,ds\right)^4
\end{split}
\eex 
for some $M\in\mathcal F$. This completes the proof. \endpf\\

\noindent{\bf Proof of Theorem \ref{le:Phi}}. It is readily seen that the conclusion
follows from (\ref{rho u_tt estimates 3}), (\ref{H^1 of u:2}), (\ref{w^1,q of
rho:2}), (\ref{exss8.6}) and (\ref{H^3 of d}). \qed

\section {Proof of Theorem \ref{App-thm:local}}
\setcounter{equation}{0} \setcounter{theorem}{0}
\subsection{$W^{2,p}$-estimate}
In this subsection, we give a proof of $W^{2,p}$-estimate of the Lam\'e equation on a simply connected,
bounded, smooth domain with the Navier-slip boundary condition, which is 
needed in our proof of Theorem \ref{App-thm:local}. We believe that such an estimate
may have its own interest.

\begin{lemma} \label{2p-est-lemma} For any simply connected, smooth bounded domain $\Omega\subset\mathbb R^3$,
$1<p<+\infty$, and $f\in L^p(\Omega,\mathbb R^3)$, If $u\in H^1\cap H^2(\Omega,\mathbb R^3)$ is 
a weak solution of 
\begin{eqnarray}\label{lame}
\mathcal L u &=& f \ {\rm{in}}\ \Omega, \nonumber\\
u\cdot\nu = (\nabla\times u)\times \nu &=& 0 \ {\rm{on}}\ \partial\Omega. 
\end{eqnarray}
Then $u\in W^{2,p}(\Omega)$, and there exists $C>0$ depending on $p, \Omega$, and $\mathcal L$ such that
\begin{equation}\label{2p-estimate}
\Big\|\nabla^2 u\Big\|_{L^p}\le C\Big[\|f\|_{L^p}+\|\nabla u\|_{L^2}\Big].
\end{equation}

\end{lemma}
\pf By the duality argument, we may assume $1<p\le 2$. Since $u\cdot\nu=0$ on $\partial\Omega$, it follows
from Bourguignon-Brezis \cite{BB} that 
\begin{equation}\label{2p-1}
\|\nabla^2 u\|_{L^p}\les \|\nabla({\rm{div}} \ u)\|_{L^p}
+\|\nabla({\rm{curl}}\ u)\|_{L^p}+\|\nabla u\|_{L^p}.
\end{equation}
Also, since $\Omega$ is simply connected and $(\nabla\times u)\times\nu=0$ on $\partial\Omega$, 
it follows from Wahl \cite{Wolf} that
\begin{eqnarray}\label{2p-2}
\|\nabla({\rm{curl}}\ u)\|_{L^p}
&\leq& C\|\nabla\times{\rm{curl}}\ u\|_{L^p}+\|\nabla\cdot({\rm{curl}}\ u)\|_{L^p}
=C\|\nabla\times ({\rm{curl}}\ u)\|_{L^p}\nonumber\\
&\lesssim& \frac{1}{\mu}\Big[\|\mathcal L u\|_{L^p}+(2\mu+\lambda)\|\nabla({\rm{div}} u)\|_{L^p}\Big]\nonumber\\
&\lesssim& \|\nabla({\rm{div}} u)\|_{L^p}+\|f\|_{L^p}.
\end{eqnarray}
Now we estimate $\|\nabla({\rm{div}} \ u)\|_{L^p}$ by the duality argument: for $p'=\frac{p}{p-1}$,
$$
\|\nabla({\rm{div}} u)\|_{L^p}
\le C\sup\Big\{\int \nabla({\rm{div}} u)\cdot g\,dx: \ g\in C^\infty(\overline\Omega, \mathbb R^3),
\ \|g\|_{L^{p'}}=1\Big\}.$$

For any $g\in C^\infty(\overline\Omega,\mathbb R^3)$, with $\|g\|_{L^{p'}}=1$, by the Helmholtz's
decomposition Theorem (see Fujiwara-Morimoto \cite{FM} and Solonnikov \cite{S}), there exist
$G\in C^\infty(\overline\Omega) \cap W^{1,p'}(\Omega)$ and $H\in C^\infty(\overline\Omega)\cap L^{p'}(\Omega,\mathbb R^3)$ such that
\begin{eqnarray*}
g&=&\nabla G+H, \ {\rm{div}} H=0 \ {\rm{in}}\ \Omega, \\
\frac{\partial G}{\partial\nu} &=& g\cdot\nu \ {\rm{on}}\ \partial\Omega,\\
\|G\|_{W^{1,p'}} &+&\|H\|_{L^{p'}}\le C\|g\|_{L^{p'}}=C.
\end{eqnarray*}
Thus we have 
$$\int \nabla({\rm{div}}\ u)\cdot H\,dx=0$$
so that
\begin{eqnarray*}
\int \nabla({\rm{div}}\ u)\cdot g \,dx&=& \int \nabla({\rm{div}}\ u)\cdot(\nabla G+H)\,dx
=\int \nabla({\rm{div}}\ u)\cdot\nabla G\,dx\\
&=&\int (\nabla({\rm{div}}\ u)-\frac{1}{2\mu+\lambda} f)\cdot\nabla G \,dx+\frac{1}{2\mu+\lambda}
\int f\cdot\nabla G\,dx\\
&=&\frac{\mu}{2\mu+\lambda}\int \nabla\times({\rm{curl}}\ u)\cdot\nabla G\,dx
+\frac{1}{2\mu+\lambda}
\int f\cdot\nabla G\,dx\\
&=&\frac{1}{2\mu+\lambda}
\int f\cdot\nabla G\,dx,
\end{eqnarray*}
where we have used
$$\int \nabla\times ({\rm{curl}}\ u)\cdot\nabla G=0,$$
since ${\rm{div}}(\nabla\times ({\rm{curl}}\ u))=0$ in $\Omega$ and 
$({\rm{curl}}\ u)\times \nu =0$ on $\partial\Omega$. 
The above inequality implies
$$\Big|\int \nabla({\rm{div}}\ u)\cdot g\,dx\Big|
\lesssim \|f\|_{L^p}\|\nabla G\|_{L^{p'}}\le C\|f\|_{L^p}.$$
Taking supremum over all such $g$'s, we obtain 
$$\|\nabla({\rm{div}}\ u)\|_{L^p}\leq C\|f\|_{L^p}.$$
It is clear that  this, with the help of (\ref{2p-1}) and (\ref{2p-2}), 
implies (\ref{2p-estimate}). \qed

\subsection{Existence}
In this subsection, we will first consider that $\Omega\subset\mathbb R^3$ is a bounded domain,  and then 
employ the Galerkin's method to obtain a sequence of approximate solutions to (\ref{clc-1})-(\ref{clc-3}) under (\ref{clcinitial}) and (\ref{clcboundary2}) or (\ref{clcboundary3}) that enjoy {\it a priori} estimates
obtained in \S2, which will converge to a strong solution to (\ref{clc-1})-(\ref{clc-3}). 
The existence of strong solutions for the Cauchy problem on $\mathbb R^3$ follows in a standard way
from a priori estimates by the {\it domain exhaustion} technique, which will be sketched at the end
of this subsection.

To implement the Galerkin's method, we take the function space $X$ to be either \\
(i) for the Dirichlet boundary condition (\ref{clcboundary2}), $X:=H^1_0\cap H^2(\Omega,\mathbb R^3)$ 
and and its finite dimensional subspaces as
$$
X^m:={\rm{span}}\left\{\phi^1,\cdots,\phi^m\right\}, \ m\ge 1,$$
where $\{\phi^m\}\subset X$  is an orthonormal base of $H^1(\Omega)$,
formed by the  set of eigenfunction of the Lam\'e operator under the boundary condition 
$u=0$ on $\partial\Omega$; or\\
(ii) for the Navier-slip boundary condition (\ref{clcboundary3}), 
$$X:=\left\{u\in H^2(\Omega,\mathbb R^3): \ u\cdot\nu=(\nabla\times u)\times\nu=0 \ {\rm{on}}\ \partial\Omega\right\},$$
and its finite dimensional subspaces as
$$X^m:={\rm{span}}\left\{\phi^1,\cdots,\phi^m\right\},\ m\ge 1,$$
where $\{\phi^m\}\subset X$ is an orthonormal base of $H^1(\Omega)$,  formed by 
the set of eigenfunction of the Lam\'e operator under the Navier-slip boundary condition $\ u\cdot\nu=(\nabla\times u)\times\nu=0 \ {\rm{on}}\ \partial\Omega$. 
By the $W^{2,p}$-estimate of Lam\'e equation under (\ref{clcboundary2}) or (\ref{clcboundary3}) (see Lemma
3.1), we see that $\{\phi^m\}\subset W^{2,p}(\Omega)$ for any $1<p<+\infty$.

Now we outline the Galerkin's scheme into several steps.\\

\noindent{\it Step 1} (modification of initial data).  For $\delta>0$, let $\rho_0^\delta=\rho_0+\delta$,
$d_0^\delta=d_0$, and $u_0^\delta\in X$ be the unique solution of
\begin{eqnarray}\label{mod_initial}
\mathcal L u_0^\delta-\nabla(P(\rho_0^\delta))-\Delta d_0\cdot\nabla d_0&=&\sqrt{\rho_0^\delta}\ g
\ \ {\rm{in}}\ \Omega, \\
u_0^\delta=0; \ {\rm{or}}\ u_0^\delta\cdot\nu = (\nabla\times u_0^\delta)\times\nu&=&0 
\ \ \ \ \ \ \ {\rm{on}}\ \partial\Omega. \label{mod_bdry}
\end{eqnarray}
By the $W^{2,2}$-estimate of Lam\'e equation, it is not hard to show that
$$\lim_{\delta\downarrow 0^+}\left\|u_0^\delta-u_0\right\|_{X}=0.$$

\noindent{\it Step 2} ($m$th approximate solutions). Fix $\delta>0$ and $3<q\le 6$.
For $m\ge 1$  and some $0<T=T(m)<+\infty$
to be determined below,  we let
$$u^m_0=\sum_{k=1}^m (u_0^\delta, \phi_k)\phi_k$$
and look for the triple 
\bex
\begin{cases} \rho^m\in C([0,T];  W^{1,q}\cap H^1)\\
u^m(x,t)=\sum\limits_{k=1}^m u^m_k(t)\phi_k(x)\in C([0, T];  W^{2,q}\cap H^2)\\
d^m\in C([0,T];  H^3(\Omega, S^2))
\end{cases}
\eex
solution of the following problem
\begin{equation}\label{approx_eqn}
\begin{cases}\rho^m_t+\nabla\cdot(\rho^m u^m)=0,\\
\left(\rho^m u^m_t,\ \phi_k\right)+\mu(\nabla \times u^m, \ \nabla\phi_k)+(2\mu+\lambda)(\nabla\cdot u^m,\ \nabla\phi_k) \\
 =-(\rho^m u^m\cdot\nabla u^m,\ \phi_k)
-(\nabla(P(\rho^m)),\ \phi_k)-(\Delta d^m\cdot\nabla d^m, \ \phi_k) \ (1\le k\le m),\\
d^m_t+u^m\cdot\nabla d^m =\Delta d^m +|\nabla d^m|^2d^m,\\
(\rho^m, \ u^m,\  d^m)\Big|_{t=0} =(\rho_0^\delta,\ u^m_0,\ d_0),\\
(u^m, \ \frac{\partial d^m}{\partial\nu})\Big|_{\partial\Omega\times [0,T]}=0,  
\ {\rm{or}}\ (u^m\cdot\nu, \ (\nabla\times u^m)\times\nu, \ \frac{\partial d^m}{\partial\nu})
\Big|_{\partial\Omega\times [0,T]}=0.
\end{cases}
\end{equation}
The existence of a solution $(\rho^m, u^m, d^m)$ to (\ref{approx_eqn}) 
over $\Omega\times [0,T(m)]$ for some $T(m)>0$ can be obtained by the fixed point theorem, similar 
to that on the compressible Navier-Stokes equation by Padula \cite{P} (see also \cite{CK}). Here we only sketch the argument. First,  observe that for any given $0<T<+\infty$ and $u^m\in C([0, T]; W^{2,q}\cap H^2)$, it is standard to show that there exist \\
(1) a solution $\rho^m\in C([0, T]; W^{1,q}\cap H^1)$  of (\ref{approx_eqn})$_1$ along with 
$\rho^m\Big|_{t=0}=\rho_0^\delta$.\\
(2) $0<t_m\le T$, depending on $u^m$ and $\|d_0\|_{H^3}$, and a solution
$d^m\in C([0, t_m], H^3(\Omega,S^2))$ of (\ref{approx_eqn})$_3$ along with
$d^m\Big|_{t=0}=d_0$ and $\frac{\partial d^m}{\partial\nu}\Big|_{\partial\Omega\times [0,t_m]}=0$. 

It is well-known (cf. \cite{P} \cite{CK} or Lemma 2.5 in \S2) that 
\beq\label{lowerbound}
\rho^m(x,t)\ge \delta \exp\left(-\int_0^t \|\nabla u^m\|_{L^\infty}\,ds\right)>0, \ (x,t)\in Q_T.
\eeq
The coefficients $u^m_k(t)$  can be determined by the following system of $m$ 
first order ordinary differential equations: $1\le k\le m$, 
\beq\label{ode}
\sum_{i=1}^m (\rho^m\phi_i, \phi_k) \dot{u}^m_i=F_k\left(u^m_l(t), \int_0^t u^m_l\,ds, t\right);
\ u^m_k(0)=(u_0^\delta,\phi_k),
\eeq
where $F_k$ denotes the right hand side of (\ref{approx_eqn})$_2$. Since $\rho^m$ is strictly positive,
the determinant of the $m\times m$ matrix $\left(\rho^m\phi_i, \phi_k\right)_{1\le i, k\le m}$ is positive. Hence we can reduce (\ref{ode})
into 
\beq{}\label{ode1}
\dot{u}^m_k=G_k(u^m_l, b^m_l, t), \  \dot{b}^m_k=u^m_k;
\ u^m_k(0) =(u_0^\delta,\phi_k), \  b^m_k(0)=0,
\eeq
where $G_k$ is a regular function of $u^m_l, b^m_l$. Therefore, by the standard existence theory of
ordinary differential equations, we conclude that there exists a $0<T_m\le t_m$ and a solution
$u^m_k(t)$ to (\ref{ode}), which in turn implies the existence of solutions $\rho^m, d^m$ 
of (\ref{approx_eqn})$_1$ and (\ref{approx_eqn})$_3$ on the same time interval.

\medskip
\noindent{\it Step 3} (a priori estimates). We will show that there exist  $0<T_0<+\infty$ and $C>0$, depending only on the norms given by
the regularity conditions on $P$ and the initial data $\rho_0, u_0$, and $d_0$, but independent of the parameters $\delta, m$, and the size of
the domain $\Omega$, such that there exists $M\in\mathcal F$ so that for any $m\ge 1$,
$(\phi^m, u^m, d^m)$ satisfies:
\begin{equation}\label{apriori-est}
\Phi^m(t)\le \exp \left[C\mathcal M(\rho_0^\delta, u_0^\delta, d_0^\delta)+
C\int_0^t M(\Phi^m(s))\,ds\right],
\ 0<t\le T_0,
\end{equation}
where $\Phi^m(t)$ is defined by (\ref{phi}) with $(\rho, u, d)$ replaced by $(\rho^m, u^m, d^m)$ and $\mathcal M(\rho_0^\delta, u_0^\delta, d_0^\delta)$ is defined by (\ref{rho-u-d})  with
$(\rho_0,u_0,d_0)$ replaced by $(\rho_0^\delta, u_0^\delta, d_0^\delta)$. 

Since the argument to obtain (\ref{apriori-est}) is almost identical to proof of Theorem \ref{le:Phi}, we only birefly outline it here:

First, it is easy to see (\ref{approx_eqn})$_2$ holds with $\phi_k$ replaced by $u^m$. By multiplying (\ref{approx_eqn})$_3$
by $(\Delta d^m+|\nabla d^m|^2 d^m)$ and integrating over $\Omega$ and adding these two resulting equations, we can show
that there is a $M\in\mathcal F$ such that the energy inequality (\ref{Energy-identity}) holds with $(\rho, u, d)$, $M$, and $\Phi$ replaced by $(\rho^m, \ u^m,\ d^m)$, $M$, and $\Phi^m$. 

Second, since (\ref{approx_eqn})$_2$ implies
\beq\label{um-eqn}
\mathcal L u^m =\mathbb P_m\Big(\rho^m \dot{u}^m+\nabla(P(\rho^m))+\nabla d^m\cdot\Delta d^m\Big),
\eeq
where $\mathbb P_m(u)=\sum_{i=1}^m (u,\phi_k)\phi_k: X\to X^m$ is the orthogonal projection map,
we can check that the same argument as Lemma \ref{le:2.2} yields that exists $M\in\mathcal F$ so that
\beq\label{H2-um}
\|\nabla u^m\|_{H^1}^2\le M(\Phi^m(t)), \ 0\le t\le T_m
\eeq

Third, by differentiating (\ref{um-eqn}) w.r.t. $t$,  multiplying the resulting equation with $u^m_t$, integrating over $\Omega$, and
repeating the proof of Lemma \ref{le:nabla u_t}, we obtain that there exists $M\in\mathcal F$ such that
for any $m\ge 1$, 
\beq\label{dotut-est}
\int \rho^m |u^m_t|^2+\int_0^t\int_\Omega |\nabla u^m_t|^2
\le C\left[\mathcal M(\rho_0^\delta, u_0^m, d_0^\delta)+\int_0^t M(\Phi^m(s))\,ds\right].
\eeq
Fourth, similar to the proof of Lemma \ref{le:2.4} and Lemma \ref{nabla rho est}, we have
that there exists $M\in\mathcal F$ such that for all $m\ge 1$, 
\beq{}\label{nablaum-est}
\|\nabla u^m\|_{L^2}^2\le C\left[\mathcal M(\rho_0^\delta, u^m_0, d_0^\delta)+\int_0^t 
M(\Phi^m(s))\,ds\right],
\eeq
and
\beq\label{nabla rho delta est}
\|\rho^m\|_{H^1\cap W^{1,q}}\le C\exp \left\{C\left[\mathcal M(\rho_0^\delta, u^m_0,d_0^\delta)+
\int_0^t M(\Phi^m(s))\,ds\right]\right\}.
\eeq
Fifth, by differentiating (\ref{approx_eqn})$_3$ w.r.t. $x$ and mutiplying by $\nabla d^m_t$ (and $\nabla\Delta d^m_t$ respectively) and integrating over $\Omega$,
we can use the same argument as Lemma \ref{le:H^2 of d} and Lemma \ref{le:H^3 of d} to show that
there exists $M\in \mathcal F$ such that for all $m\ge 1$, 
\beq\label{H^2 of dm}
\|\nabla^2 d^m\|_{L^2}^2+\int_0^t \|\nabla d^m_t\|_{L^2}^2\,ds
\le C[1+\int_0^t M(\Phi^m(s))\,ds],
\eeq
\beq\label{H^3 of dm}
\|\nabla^3 d^m\|_{L^2}^2+\int_0^t\|\nabla^2 d^m_t\|_{L^2}^2\,ds
\le \left(C\mathcal M(\rho_0^\delta, u^m_0,d_0^\delta)+\int_0^t M(\Phi^m(s))\,ds\right)^4.
\eeq
It is readily seen that combining all these estimates together yields
(\ref{apriori-est}) with $T_0$ replaced by $T_m$ and $u_0^\delta$ replaced by $u_0^m$.  \\

\noindent{\it Step 4} (convergence and solution).  By the definition of $u^\delta_0$, $\mathcal M$ given by
(\ref{phi}),  and the condition (\ref{first3.1}),  we have
$$\mathcal M(\rho_0^\delta, u_0^\delta, d_0^\delta)
=1+\|g\|_{L^2},$$
and
$$\Big|\mathcal M(\rho_0^\delta, u^m_0, d_0^\delta)
-\mathcal M(\rho_0^\delta, u_0^\delta, d_0^\delta)\Big|
\le \frac{C}{\delta}\Big\|u_0^m-u_0^\delta\Big\|_{H^2}\rightarrow 0, 
\ {\rm{as}}\ m\rightarrow\infty. $$
Thus there exists $N=N(\delta)>0$ such that
\beq\label{mbound}
\mathcal M(\rho_0^\delta, u_0^m, d_0^\delta)\le 2+\|g\|_{L^2},
\ \forall m\ge N.
\eeq
It follows from (\ref{mbound}), (\ref{apriori-est}), and Gronwall's inequality (see, for example, 
\cite{CCK} page 263 or \cite{Simon} Lemma 6)
that there exists a small $T_0>0$, independent of $\delta$ and $m$, such that
\beq\label{phi_m-bound}
\sup_{0\le t\le T_0} \Phi^m(t)
\le C\exp (C\|g\|_{L^2}), \ \forall m\ge M.
\eeq 
By virtue of (\ref{phi_m-bound}), we obtain that for any $m\ge M$,
\begin{eqnarray}\label{strong_sol-est}
&&\sup_{0\le t\le T_0}
\Big(\|\sqrt{\rho^m} u^m_t\|_{L^2}^2+\|\rho^m\|_{W^{1,q}\cap H^1}^2
+\|\nabla u^m\|_{H^1}^2+\|d^m_t\|_{H^1}^2+\|\nabla d^m\|_{H^2}^2\Big)\nonumber\\
&+&\int_0^{T_0}\Big(\|u^m\|_{D^{2,q}}^2+\|\nabla u^m_t\|_{L^2}^2
+\|\nabla^4 d^m\|_{L^2}^2+\|\nabla^2 d_t^m\|_{L^2}^2\Big)
\le  C\exp (C\|g\|_{L^2}^2).
\end{eqnarray}
Based on the estimate (\ref{strong_sol-est}), we can deduce that after taking subsequences, 
there exists $(\rho^\delta, u^\delta, d^\delta)$ such that
$$\rho^m\rightharpoonup \rho^\delta 
\ {\rm{weak}}^* \ {\rm{in}}\ L^\infty(0, T_0; W^{1,q}\cap H^1), 
\ u^m\rightharpoonup u^\delta 
\ {\rm{weak}}^* \ {\rm{in}}\ L^\infty(0, T_0; D^1\cap D^2),$$
$$
\ u^m\rightharpoonup u^\delta 
\ {\rm{weak}}\ {\rm{in}}\ L^2(0, T_0; D^{2,q}),
\ u^m_t\rightharpoonup u^\delta_t
\ {\rm{weak}}\ {\rm{in}}\ L^2(0, T_0; D^1),
$$
$$
\sqrt{\rho^m}u^m_t\rightharpoonup \sqrt{\rho^\delta} u^\delta_t 
\ {\rm{weak}}^*\ {\rm{in}}\ L^\infty(0,T; L^2),$$
$$
d^m\rightharpoonup d^\delta 
\ {\rm{weak}}^* \ {\rm{in}}\ L^\infty(0, T_0; D^1\cap D^3) \ {\rm{and}}\ 
L^2(0,T_0; D^4),$$
$$ \ d_t^m\rightharpoonup d^\delta_t \ {\rm{in}}\ L^2(0, T; H^2)\ 
\mathrm{and\ weak}^*\ \mathrm{in}\ L^\infty(0,T; H^1).
$$
By the lower semicontinuity, (\ref{strong_sol-est}) implies that for $0\le t\le T_0$,
$(\rho^\delta, u^\delta, d^\delta)$ satisfies
\begin{eqnarray}\label{strong_sol-est1}
&&\sup_{0\le t\le T_0}
\Big(\|\sqrt{\rho^\delta} u^\delta_t\|_{L^2}^2+\|\rho^\delta\|_{W^{1,q}\cap H^1}^2
+\|\nabla u^\delta\|_{H^1}^2+\|d^\delta_t\|_{H^1}^2+\|\nabla d^\delta\|_{H^2}^2\Big)\nonumber\\
&+&\int_0^{T_0}\Big(\|u^\delta\|_{D^{2,q}}^2+\|\nabla u^\delta_t\|_{L^2}^2
+\|\nabla^4 d^\delta\|_{L^2}^2+\|\nabla^2 d_t^\delta\|_{L^2}^2\Big)
\le  C\exp (C\|g\|_{L^2}^2).
\end{eqnarray}
Furthermore, it is straightforward to check that $(\rho^\delta, u^\delta, d^\delta)$
is a strong solution in $[0, T_0]$ of (\ref{clc-1})-(\ref{clc-3}) under the initial condition
$(\rho^\delta,u^\delta,d^\delta)\Big|_{t=0}=(\rho^\delta_0, u^\delta_0, d^\delta_0)$
and the boundary condition (\ref{clcboundary2}) or (\ref{clcboundary3}). Since
$T_0>0$ is independent of $\delta$,  $(\rho^\delta, \ u^\delta,\ d^\delta)$ satisfies
(\ref{strong_sol-est1}), 
$\rho_0^\delta\rightarrow \rho_0$ in $W^{1,q}\cap H^1$,
$u_0^\delta\rightarrow u_0$ in $D^1\cap D^2$, and $d_0^\delta=d_0$, the same
limiting process as above would imply that after taking a subsequence $\delta\downarrow 0$,
$(\rho^\delta,u^\delta, d^\delta)$ converges (weakly in the corresponding spaces)
to a strong solution $(\rho, u,d)$  of (\ref{clc-1})-(\ref{clc-3}) on $\Omega\times [0,T_0]$
along with (\ref{clcinitial}) and
(\ref{clcboundary2}) or (\ref{clcboundary3}). 

For the Cauchy problem on $\mathbb R^3$,  we proceed as follows.
For $R\uparrow\infty$,  it is standard (cf. \cite{LW1}) that there exists $d_0^R\in H^3(\mathbb R^3, S^2)$
such that $d_0^R \equiv n_0$ outside $B_{\frac{R}2}$ for some constant $n_0\in S^2$ and
\beq\label{d_0-approx}
\lim_{R\uparrow\infty}\Big\|\nabla d_0^R-\nabla d_0\Big\|_{H^2(\mathbb R^3)}=0.
\eeq
Now we let $u_0^R\in H^1_0(B_R)\cap H^2(B_R)$ be the unique solution of
\beq\label{lame_approx}
\mathcal L u_0^R-\nabla (P(\rho_0))-\Delta d_0^R\cdot\nabla d_0^R
=\sqrt{\rho_0} \ g \ {\rm{on}}\ B_R,\ u_0^R\Big|_{\partial B_R}=0,
\eeq
where $g\in L^2(\mathbb R^3)$ is given by (\ref{first3.1}).  Extending $u_0^R$ to $\mathbb R^3$ by
letting it be zero outside $B_R$.  Then it is not hard to show that for any compact subset $K\subset\mathbb R^3$,
\beq\label{u_0-approx}
\lim_{R\uparrow \infty} \Big\|\nabla u_0^R-\nabla u_0\Big\|_{H^1(K)}=0.
\eeq

By the above existence, we know that there exists $T_0>0$, independent of $R$,
and  a strong solution  $(\rho^R,\ u^R,\ d^R)$ of (\ref{clc-1})-(\ref{clc-3}) 
on $B_R\times [0,T_0]$ of (\ref{clc-1})-(\ref{clc-3}),
under the initial and boundary condition:
\beq\label{ivp-R}
(\rho^R,\ u^R,\ d^R)\Big|_{B_R\times \{t=0\}}=(\rho_0,\ u_0^R,\ d_0^R);
\ (u^R,\ \frac{\partial d^R}{\partial R})\Big|_{\partial B_R\times [0,T_0]}=0.
\eeq
Furthermore, $(\rho^R, u^R, d^R)$ satisfies the estimate:
\begin{eqnarray}\label{strong_sol-est2}
&&\sup_{0\le t\le T_0}
\Big(\|\sqrt{\rho^R} u^R_t\|_{L^2}^2+\|\rho^R\|_{W^{1,q}\cap H^1}^2
+\|\nabla u^R\|_{H^1}^2+\|d^R_t\|_{H^1}^2+\|\nabla d^R\|_{H^2}^2\Big)\nonumber\\
&+&\int_0^{T_0}\Big(\|u^R\|_{D^{2,q}}^2+\|\nabla u^R_t\|_{L^2}^2
+\|\nabla^4 d^R\|_{L^2}^2+\|\nabla^2 d_t^R\|_{L^2}^2\Big)
\le  C\exp (C\|g\|_{L^2}^2),
\end{eqnarray}
with $C>0$ independent of $R$.  It is readily seen that (\ref{strong_sol-est2}), (\ref{d_0-approx}),
and (\ref{u_0-approx}) imply that after taking a subsequence, we may assume that
$(\rho^R, u^R, d^R)$ locally converges (weakly in the corresponding spaces)
to a strong solution $(\rho, u, d)$ of (\ref{clc-1})-(\ref{clc-3})  on
$\mathbb R^3\times [0,T_0]$ under the initial condition
(\ref{clcinitial}) and the boundary condition (\ref{clcboundary1}).  This completes the proof of Theorem
\ref{App-thm:local}. \qed

\subsection{Uniqueness}
In this subsection, we will show the uniqueness of the local strong solutions obtained in 
Theorem \ref{App-thm:local}.

Let $(\rho_i,u_i,d_i)$ ($i=1,2$) be two strong solutions on $\Omega\times (0,T]$  of
(\ref{clc-1})-(\ref{clc-3}) with (\ref{clcinitial}) and either (\ref{clcboundary1}),
or (\ref{clcboundary2}), or (\ref{clcboundary3}). Set
$\overline{\rho}=\rho_2-\rho_1,\overline{u}=u_2-u_1,\overline{d}=d_2-d_1$. Then
we have  \beq
  \label{3.45}\begin{cases}
          \overline{\rho}_t+(u_1\cdot
\nabla)\overline{\rho}+\overline{u}\cdot\nabla\rho_2+\overline{\rho}\mathrm{div}u_2+\rho_1\mathrm{div}\overline{u}=0,\\[2mm]
          \rho_1\overline{u}_t+\rho_1u_1\cdot\nabla \overline{u}+\nabla \left(P(\rho_2)-P(\rho_1)\right)\\
=\mathcal L\overline u-\overline{\rho}(u_{2t}+u_2\cdot\nabla u_2)
  -\rho_1\overline{u}\cdot\nabla u_2-\Delta\overline{d}\cdot\nabla d_2-\Delta d_1\cdot\nabla\overline{d},\\[2mm]
        \overline{d}_t-\triangle \overline{d}=\nabla \overline{d}\cdot(\nabla d_2+\nabla d_1)d_1+|\nabla
        d_2|^2\overline{d}
        -\overline{u}\cdot\nabla d_2-u_1\cdot\nabla\overline{d},\\[1mm]
\end{cases}
\eeq with the initial condition:
 \bex
 (\overline{\rho},\overline{u},\overline{d})|_{t=0}=0,\ x\in\overline{\Omega},\eex
and the boundary condition: \bex
        \big(\overline{u},\frac{\partial\overline{d}}{\partial
        \nu}\big)\big|_{\partial\Omega}=0,
\ {\rm{or}}\ 
        \big(\overline{u}\cdot\nu, (\nabla\times\overline{u})\times\nu,\frac{\partial\overline{d}}{\partial
        \nu}\big)\big|_{\partial\Omega}=0.
\eex Multiplying $(\ref{3.45})_2$ by $\overline{u}$, integrating
over $\Omega$, and using integration by parts, we have \bex
&&\frac{1}{2}\frac{d}{dt}\int\rho_1|\overline{u}|^2\,dx+
\int\left((2\mu+\lambda)|\mathrm{div}\overline{u}|^2+\mu|\nabla\times\overline{u}|^2\right)\,dx\\
&=&
-\int\overline{\rho}(u_{2t}+u_2\cdot\nabla u_2)
\cdot\overline{u}\,dx -\int\rho_1\overline{u}\cdot\nabla
u_2\cdot\overline{u}\,dx+\int\left(P(\rho_2)-P(\rho_1)\right)\mathrm{div}\overline{u}\,dx\\
&+&\int\left(\nabla\overline{d}\cdot\nabla\nabla
d_2\cdot\overline{u}+\nabla\overline{d}\cdot\nabla
d_2\cdot\nabla\overline{u}\right)\,dx-\int\Delta
d_1\cdot\nabla\overline{d}\cdot\overline{u}\,dx. \eex
Observe that
$$|P(\rho_2)-P(\rho_1)|\le B_P\Big(\|\rho_1\|_{L^\infty}+\|\rho_2\|_{L^\infty}\Big)|\overline\rho|
\le C|\overline\rho|.$$
Hence, by H\"older's inequality and Cauchy's inequality, we have
\bex
&&\frac{1}{2}\frac{d}{dt}\int\rho_1|\overline{u}|^2dx+
\int\left((2\mu+\lambda)|\mathrm{div}\overline{u}|^2+\mu|\nabla\times\overline{u}|^2\right)dx\\&\lesssim&\|\overline{\rho}\|_{L^\frac{3}{2}}
\|u_{2t}+u_2\cdot\nabla u_2\|_{L^6}\|\overline{u}\|_{L^6}+\|\nabla
u_2\|_{L^\infty}\int_\Omega\rho_1|\overline{u}|^2dx+\|\overline{\rho}\|_{L^2}\|\mathrm{div}\overline{u}\|_{L^2}\\
&+&
\|\nabla\overline{d}\|_{L^2}\|\nabla^2
d_2\|_{L^3}\|\overline{u}\|_{L^6}+\|\nabla\overline{d}\|_{L^2}\|\nabla
d_2\|_{L^\infty}\|\nabla\overline{u}\|_{L^2}+\|\Delta
d_1\|_{L^3}\|\nabla\overline{d}\|_{L^2}\|\overline{u}\|_{L^6}\\
&\lesssim&\|\overline{\rho}\|_{L^\frac{3}{2}}
\|u_{2t}+u_2\cdot\nabla u_2\|_{L^6}\|\nabla\overline{u}\|_{L^2}
+\|\nabla
u_2\|_{W^{1,q}}\int\rho_1|\overline{u}|^2dx\\
&+&\|\overline{\rho}\|_{L^2}\|\mathrm{div}\overline{u}\|_{L^2}
+\|\nabla\overline{d}\|_{L^2}\|\nabla\overline{u}\|_{L^2}\\&\le&
\epsilon\int |\nabla\overline u|^2\,dx+C[\|\overline{\rho}\|_{L^\frac{3}{2}}^2
\|u_{2t}+u_2\cdot\nabla u_2\|_{L^6}^2\\
&+&\|\nabla
u_2\|_{W^{1,q}}\int\rho_1|\overline{u}|^2dx+\|\overline{\rho}\|_{L^2}^2
+\|\nabla\overline{d}\|_{L^2}^2]. \eex 
Thus, by choosing $\epsilon$ sufficiently small, we have
\beq\label{3.46}\begin{split}
&\frac{d}{dt}\int\rho_1|\overline{u}|^2\,dx+
\int |\nabla\overline u|^2\,dx\\
\le& C[\|\overline{\rho}\|_{L^\frac{3}{2}}^2 \|u_{2t}+u_2\cdot\nabla
u_2\|_{L^6}^2+\|\nabla
u_2\|_{W^{1,q}}\int\rho_1|\overline{u}|^2dx+\|\overline{\rho}\|_{L^2}^2
+\|\nabla\overline{d}\|_{L^2}^2].\end{split}\eeq 
Multiplying
$(\ref{3.45})_1$ by $2\overline{\rho}$, integrating over $\Omega$,
and using integration by parts, we have
\beq\label{3.47}\begin{split}
\frac{d}{dt}\int|\overline{\rho}|^2\,dx\lesssim&\int|\overline{\rho}\
\overline{u}\cdot\nabla\rho_2|\,dx
+\int|\overline{\rho}|^2(|\mathrm{div}u_1|+|\mathrm{div}u_2|)\,dx+\int|\overline{\rho}\
\rho_1\mathrm{div}\overline{u}|\,dx\\
\lesssim&\|\overline{\rho}\|_{L^2}
\|\nabla\rho_2\|_{L^3}\|\overline{u}\|_{L^6}+(\|\mathrm{div}
u_1\|_{L^\infty}+\|\mathrm{div}u_2\|_{L^\infty})\int|\overline{\rho}|^2\,dx+
\|\overline{\rho}\|_{L^2}\|\mathrm{div}\overline{u}\|_{L^2}\\
\lesssim&\|\overline{\rho}\|_{L^2}
\|\nabla\overline{u}\|_{L^2}+(\|\mathrm{div}u_1\|_{W^{1,q}}+\|\mathrm{div}u_2\|_{W^{1,q}})\int|\overline{\rho}|^2\,dx\\
\lesssim&\|\overline{\rho}\|_{L^2}(\|\mathrm{div}\overline{u}\|_{L^2}+\|\nabla\times\overline{u}\|_{L^2})+(\|\mathrm{div}u_1\|_{W^{1,q}}
+\|\mathrm{div}u_2\|_{W^{1,q}})
\int |\overline{\rho}|^2dx\\
\le&\epsilon\int |\nabla\overline u|^2\,dx+C_\epsilon\|\overline{\rho}\|_{L^2}^2+C(\|\mathrm{div}u_1\|_{W^{1,q}}+\|\mathrm{div}u_2\|_{W^{1,q}})\int|\overline{\rho}|^2\,dx,\end{split}
\eeq 
for any $\epsilon>0$.
Similarly, we have \beq\label{rho 3/2:1}
\begin{split}
&\frac{d}{dt}\int|\overline{\rho}|^\frac{3}{2}\,dx
\lesssim\int|\overline{\rho}^\frac{1}{2}\
\overline{u}\cdot\nabla\rho_2|\,dx
+\int |\overline{\rho}|^\frac{3}{2}(|\mathrm{div}u_1|+|\mathrm{div}u_2|)\,dx
+\int|\overline{\rho}^\frac{1}{2}\rho_1\mathrm{div}\overline{u}|\,dx\\
\lesssim&\|\overline{\rho}\|_{L^\frac{3}{2}}^\frac{1}{2}
\|\nabla\rho_2\|_{L^2}\|\overline{u}\|_{L^6}+(\|\mathrm{div}
u_1\|_{L^\infty}+\|\mathrm{div}u_2\|_{L^\infty})\int|\overline{\rho}|^\frac{3}{2}\,dx+
\|\overline{\rho}\|_{L^\frac{3}{2}}^\frac{1}{2}\|\mathrm{div}\overline{u}\|_{L^2}\|\rho_1\|_{L^6}\\
\lesssim&\|\overline{\rho}\|_{L^\frac{3}{2}}^\frac{1}{2}
\|\nabla\overline{u}\|_{L^2}+(\|\mathrm{div}u_1\|_{W^{1,q}}+\|\mathrm{div}u_2\|_{W^{1,q}})\int|\overline{\rho}|^\frac{3}{2}dx\\
\lesssim&\|\overline{\rho}\|_{L^\frac{3}{2}}^\frac{1}{2}\|\nabla\overline u\|_{L^2}
+(\|\mathrm{div}u_1\|_{W^{1,q}}
+\|\mathrm{div}u_2\|_{W^{1,q}})
\int|\overline{\rho}|^\frac{3}{2}dx.
\end{split}
\eeq 
Multiplying (\ref{rho 3/2:1}) by
$\|\overline{\rho}\|_{L^\frac{3}{2}}^\frac{1}{2}$, and using
Cauchy's inequality, we have \beq\label{rho 3/2}
\begin{split}
&\frac{d}{dt}\|\overline{\rho}\|_{L^\frac{3}{2}}^2
\lesssim\|\overline{\rho}\|_{L^\frac{3}{2}}\|\nabla\overline u\|_{L^2}+(\|\mathrm{div}u_1\|_{W^{1,q}}
+\|\mathrm{div}u_2\|_{W^{1,q}})
\|\overline{\rho}\|_{L^\frac{3}{2}}^2
\\ \le&\epsilon\int|\nabla\overline u|^2\,dx
+C_\epsilon\|\overline{\rho}\|_{L^\frac{3}{2}}^2+C(\|\mathrm{div}u_1\|_{W^{1,q}}+\|\mathrm{div}u_2\|_{W^{1,q}})\|\overline{\rho}\|_{L^\frac{3}{2}}^2.
\end{split}
\eeq Multiplying (\ref{3.45})$_3$ by $-\triangle \overline{d}$,
integrating over $\Omega$, and using integration by parts and
Cauchy's inequality, we have \bex\begin{split}
\frac{1}{2}\frac{d}{dt}\int |\nabla\overline{d}|^2\,dx+\int|\triangle
\overline{d}|^2\,dx \lesssim&\|\nabla \overline{d}\|_{L^2}\|\Delta
\overline{d}\|_{L^2}\|\nabla d_2+\nabla d_1\|_{L^\infty}+\|\Delta
\overline{d}\|_{L^2}\|\nabla d_2\|_{L^6}^2\|\overline{d}\|_{L^6}\\
&+\|\Delta \overline{d}\|_{L^2}\|\overline{u}\|_{L^6}\|\nabla d_2\|_{L^3}+\|\Delta \overline{d}\|_{L^2}\|u_1\|_{L^\infty}\|\nabla\overline{d}\|_{L^2}
        \\ 
\lesssim&\|\nabla \overline{d}\|_{L^2}\|\Delta \overline{d}\|_{L^2}\|\nabla d_2+\nabla d_1\|_{H^2}+\|\Delta \overline{d}\|_{L^2}\|\nabla
        d_2\|_{H^1}^2\|\nabla\overline{d}\|_{L^2}
        \\&+\|\Delta \overline{d}\|_{L^2}\|\nabla\overline{u}\|_{L^2}\|\nabla d_2\|_{L^2}^\frac{1}{2}\|\nabla d_2\|_{L^6}^\frac{1}{2}+\|\Delta \overline{d}\|_{L^2}\|\nabla u_1\|_{H^1}\|\nabla\overline{d}\|_{L^2} \\ \lesssim&\|\Delta \overline{d}\|_{L^2}\|\nabla\overline{d}\|_{L^2}
        +\|\Delta \overline{d}\|_{L^2}\|\nabla\overline{u}\|_{L^2}\\ \le&\frac{1}{2}\|\Delta \overline{d}\|_{L^2}^2+C\|\nabla\overline{d}\|_{L^2}^2
        +C\|\nabla\overline{u}\|_{L^2}^2.
\end{split}
\eex This gives \beq\label{3.48}\begin{split}
\frac{d}{dt}\int_\Omega|\nabla\overline{d}|^2dx+\int_\Omega|\triangle
\overline{d}|^2dx \le& C\|\nabla\overline{d}\|_{L^2}^2
        +C\|\nabla\overline{u}\|_{L^2}^2.
\end{split}
\eeq 
Multiplying (\ref{3.46}) by $3C$, putting the resulting inequality,
(\ref{3.47}) and (\ref{rho 3/2}) to (\ref{3.48}), and taking
$\epsilon>0$ small enough, we have \beq\label{3.49}\begin{split}
&\frac{d}{dt}\left(3C\|\sqrt{\rho_1}\overline{u}\|_{L^2}^2+\|\overline{\rho}\|_{L^\frac{3}{2}}^2+\|\overline{\rho}\|_{L^2}^2+\|\nabla\overline{d}\|_{L^2}^2\right)
+C\int|\nabla\overline u|^2\,dx\\
\lesssim& \|\overline{\rho}\|_{L^\frac{3}{2}}^2
\|u_{2t}+u_2\cdot\nabla u_2\|_{L^6}^2+\|\nabla
u_2\|_{W^{1,q}}\int\rho_1|\overline{u}|^2dx+\|\overline{\rho}\|_{L^2}^2
+\|\nabla\overline{d}\|_{L^2}^2
\\+&(\|\mathrm{div}u_1\|_{W^{1,q}}+\|\mathrm{div}u_2\|_{W^{1,q}})\int|\overline{\rho}|^2dx+
\|\overline{\rho}\|_{L^\frac{3}{2}}^2+(\|\mathrm{div}u_1\|_{W^{1,q}}+\|\mathrm{div}u_2\|_{W^{1,q}})\|\overline{\rho}\|_{L^\frac{3}{2}}^2
\\
\lesssim& \left(\|u_{2t}+u_2\cdot\nabla u_2\|_{L^6}^2+\|\nabla
u_1\|_{W^{1,q}}+\|\nabla u_2\|_{W^{1,q}}+1\right)\\
&\cdot\left(3C\|\sqrt{\rho_1}\overline{u}\|_{L^2}^2+\|\overline{\rho}\|_{L^\frac{3}{2}}^2+\|\overline{\rho}\|_{L^2}^2+\|\nabla\overline{d}\|_{L^2}^2\right).
\end{split}\eeq
By (\ref{3.49}), Gronwall's inequality, and $(\overline{\rho_0},\ \overline{u_0},\ \overline{d_0})=0$,
we have
\beq\label{3.50}
\|\sqrt{\rho_1}\overline{u}\|_{L^2}^2+\|\overline{\rho}\|_{L^\frac{3}{2}}^2+\|\overline{\rho}\|_{L^2}^2+\|\nabla\overline{d}\|_{L^2}^2
+\int_0^t\int_\Omega|\nabla\overline{u}|^2\,dxds=0.\eeq 
This yields
\be\label{3.51}
(\overline{\rho},\ \overline{u},\ \nabla\overline{d})=0. \ee 
To see $\overline{d}=0$, observe that after substituting (\ref{3.51}) into (\ref{3.45})$_3$, we have \bex
 \overline{d}_t=|\nabla
        d_2|^2\overline{d}, \ \overline{d}|_{t=0}=0.
\eex This implies $\overline{d}=0$. 
This completes the proof. \qed

\section {Proof of Theorem \ref{umaintheorem}}
\setcounter{equation}{0} \setcounter{theorem}{0}
Let $0<T_*<\infty$ be the maximum time for the existence of strong solution $(\rho, u,d)$
to (\ref{clc-1})-(\ref{clc-3}). Namely,  $(\rho, u, d)$ is a strong
solution to (\ref{clc-1})-(\ref{clc-3}) in $\Omega\times (0, T]$ for any $0<T<T_*$, but not
a strong solution in $\Omega\times (0, T_*]$.  Suppose that (\ref{clcblpcondition}) were false,
i.e.
\begin{equation}\label{2.1}
\lim\sup\limits_{T\nearrow
T_*}\left(\|\rho\|_{L^{\infty}(0,T;L^{\infty})}+\int_0^T\|\nabla
d(t)\|^3_{L^{\infty}}\,dt\right)= M_0<\infty.
\end{equation}
The goal is to show that under the assumption (\ref{2.1}),
there is a bound $C>0$ depending only on $M_0, \rho_0, u_0, d_0$, and $T_*$ such that
\beq{}\label{uniform_est1}
\sup_{0\le t<T_*}\left[\max_{r=2, q}(\|\rho\|_{W^{1,r}}+\|\rho_t\|_{L^r})
+(\|\sqrt{\rho}u_t\|_{L^2}+\|\nabla u\|_{H^1})+(\|d_t\|_{H^1}+\|\nabla d\|_{H^2})\right]\le C,
\eeq
and
\beq\label{uniform_est2}
\int_{0}^{T_*}\left(\|u_t\|_{D^1}^2+\|u\|_{D^{2,q}}^2+\|d_t\|_{H^2}^2
+\|\nabla d\|_{H^3}^2\right)\,dt\le C.
\eeq
With (\ref{uniform_est1}) and (\ref{uniform_est2}), we can then show without much difficulty
that $T_*$ is not the maximum time, which is the desired contradiction.

The proof is based on several Lemmas.

\blm\label{uclemma2.1}{\it  Assume (\ref{2.1}), we have
\beq\label{clcblp2.2}
\int_0^{T_*}\int_{\om}|\nabla^2 d|^2\,dx\,dt\leq C. \eeq}
 \elm

\pf 
To see (\ref{clcblp2.2}), observe that (\ref{2.1}) implies  $\int_0^{T_*}\|\nabla d\|_{L^{\infty}}^2\,dt\le M_0$ so that
\bex
\int_0^{T_*}\int |\nabla d|^4\,dx\,dt&\leq& M_0
\cdot\left(\sup\ls_{0\leq t<T_*}\int |\nabla d|^2\,dx\right)\\
&\leq& M_0\Big[\int_0^{T_*} \int |P(\rho)|^2\,dxdt+\int\left(\rho_0|u_0|^2+|\nabla d_0|^2\right)\,dx\Big]
\eex
where we have used (\ref{energy_ineq11}) in the last step. Applying (\ref{energy_ineq11}) again,
this then implies
\bex
\int_0^{T_*}\int|\de d|^2\,dx\,dt&=&\int_0^{T_*}\int\left|\de d+|\nabla d|^2d\right|^2\,dx\,dt
+\int_0^{T_*}\int|\nabla d|^4\,dx\,dt\\
&\leq& (1+M_0)\Big[\int_0^{T_*}\int |P(\rho)|^2\,dxdt+\int\left(\rho_0|u_0|^2+|\nabla d_0|^2\right)\,dx
\Big].
\eex
Since
$$|P(\rho)|\le B_P(\|\rho\|_{L^\infty})|\rho|\le C|\rho|,$$
we have, by the conservation of mass and (\ref{2.1}),
\bex
\int_0^{T_*}\int |P(\rho)|^2\,dxdt
\le CT_{*} \sup_{0\le t<T_*}\|\rho\|_{L^1}\|\rho\|_{L^\infty}\le C.
\eex
Thus the standard $L^2$-estimate yields (\ref{clcblp2.2}).
\endpf\\

Following the argument by \cite{Sun-Wang-Zhang}, we let $v=\mathcal L^{-1}\nabla (P(\rho))$ be the solution of the Lam$\acute{\rm e}$ system:
\beq\label{clcblp2.6}
\begin{cases}
\mathcal L v=\nabla (P(\rho)),\\
v\big|_{\pa\om}=0,\ \mathrm{or}\ v\rightarrow0\ \mathrm{as}\ |x|\rightarrow\infty (\mathrm{when}
\ \Omega=\mathbb R^3) .
\end{cases}
\eeq
Then it follows from \cite{Sun-Wang-Zhang}  Proposition 2.1 that 
\beq \label{w^{1,q} of v}
\|\nabla v\|_{L^{q}}\le C\|P(\rho)\|_{L^q}\le CB_P(\|\rho\|_{L^\infty})\|\rho\|_{L^q}
\le C, \ 1<q\le 6,
\eeq
where we have used (\ref{2.1}) and the conservation of mass in the last step.

Denote $w=u-v$, then $w$ satisfies
\beq\label{clcblp2.7}
\begin{cases}
\rho w_t-\mathcal Lw=\rho F-\nabla d\cdot\de d,\\
 w|_{t=0}=w_0=u_0-v_0,\\ w\big|_{\pa\om}=0\ \mathrm{or}\ w\rightarrow0,\ \mathrm{as}\ |x|\rightarrow\infty,
\end{cases}
\eeq
where
\bex
\begin{split}
F=-u\cdot\nabla u-\mathcal L^{-1}\nabla(\pa_t(P(\rho)))=-u\cdot\nabla u+\mathcal L^{-1}\nabla\mbox{div
}(P(\rho) u)-\mathcal L^{-1}\nabla\big((P-P'(\rho)\rho)\mbox{div }u\big).
\end{split}
\eex
Then we have the following estimate.
\blm\label{clclemma2.2}{\it  Under the assumptions of Theorem \ref{umaintheorem},   
if $\lambda<\frac{7\mu}{9}$, then
$(\rho,u,d)$ satisfies that for any $0\leq t< T_*$,
\beq\label{clcblp2.8}
\begin{split}
\int_{\om}\left(\rho|u|^5+|\nabla w|^2+|\nabla d|^5+|\nabla ^2 d|^2\right)dx+\int_0^t\int_{\om}\left(|\nabla^3 d|^2+|\nabla^2 w|^2+|\nabla d_t|^2\right)dxds\leq C.
\end{split}
\eeq}
\elm

\pf The proof of this lemma is divided into five steps. \\
{\noindent\it Step 1. Estimates of $\displaystyle\int|\nabla w|^2\,dx$}.
Multiplying (\ref{clcblp2.7})$_1$ by $w_t$, integrating over $\om$, and using integration by parts and Cauchy's inequality, we have
\beq\label{clcblp2.10}
\begin{split}
&\frac{d}{dt}\int \left(\mu|\nabla w|^2+(\mu+\lambda)|\mbox{div }w|^2\right)dx+\int \rho|w_t|^2dx\\
\le& \|\sqrt{\rho}F\|_{L^2}^2+2\frac{d}{dt}\int (\nabla d\otimes\nabla d-\frac12{|\nabla d|^2}
\mathbb I_3):\nabla wdx+C\int |\nabla d||\nabla d_t||\nabla w|dx
=\sum\limits_{i=1}^3I_i.
\end{split}
\eeq
For $I_1$, we have
\beq\label{clcblp2.17}\begin{split}
I_1\lesssim& \|\sqrt{\rho}u\cdot\nabla u\|_{L^2}^2+\|\sqrt{\rho} \mathcal L^{-1}\nabla\mbox{div }(P(\rho) u)\|_{L^2}^2
+\|\sqrt{\rho}\mathcal L^{-1}\nabla((P(\rho)-P'(\rho)\rho)\mbox{div }u)\|_{L^2}^2\\=&
\sum\limits_{j=1}^3I_{1j}.
\end{split}\eeq
For $I_{11}$,  by H\"older's inequality, (\ref{2.1}), Sobolev inequality, interpolation inequality, and
(\ref{w^{1,q} of v}), we have
\beq\label{I_1.1-1}
\begin{split}
I_{11}\les& \|\rho^{\frac{1}{5}}u\|_{L^{5}}^2 \|\nabla u\|_{L^{\frac{10}{3}}}^2
\les\|\rho^{\frac{1}{5}}u\|_{L^{5}}^2 \|\nabla u\|_{L^{2}}^{\frac{4}{5}}\|\nabla u\|_{L^{6}}^{\frac{6}{5}}\\
\les&\|\rho^{\frac{1}{5}}u\|_{L^{5}}^2 \|\nabla u\|_{L^{2}}^{\frac{4}{5}}\|\nabla w\|_{L^{6}}^{\frac{6}{5}}+\|\rho^{\frac{1}{5}}u\|_{L^{5}}^2 \|\nabla u\|_{L^{2}}^{\frac{4}{5}}\|\nabla v\|_{L^{6}}^{\frac{6}{5}}\\
\les&\|\rho^{\frac{1}{5}}u\|_{L^{5}}^2 \|\nabla u\|_{L^{2}}^{\frac{4}{5}}\left(\|\nabla^2 w\|_{L^{2}}^{\frac{6}{5}}+\|\nabla w\|_{L^{2}}^{\frac{6}{5}}+1\right).
\end{split}\eeq
Again by \cite{Sun-Wang-Zhang} Proposition 2.1, and (\ref{clcblp2.7}), we have
\beq\label{H^2 of w}\begin{split}\|\nabla^2 w\|_{L^{2}}\les \|\sqrt{\rho}w_t\|_{L^{2}}+\|\sqrt{\rho}F\|_{L^{2}}+\|\nabla d\cdot\de d\|_{L^{2}}.\end{split}\eeq
Substituting (\ref{H^2 of w}) into (\ref{I_1.1-1}), and using Young's inequality, we obtain for any $\va>0$
\beq\label{I_1.1-3}
\begin{split}
I_{11}\le& \va(\|\sqrt{\rho}w_t\|_{L^{2}}^2+\|\sqrt{\rho}F\|_{L^{2}}^2)
\\+
&C(\|\rho^{\frac{1}{5}}u\|_{L^5}^5\|\nabla u\|_{L^{2}}^2+\|\nabla w\|_{L^2}^2+\|\nabla d\|_{L^{\infty}}^2\|\de d\|_{L^2}^2+1).
\end{split}\eeq
For $I_{12}$ and $I_{13}$, by \cite{Sun-Wang-Zhang} Proposition 2.1, (\ref{2.1}), (\ref{Energy-identity}), and
(\ref{regularity_p}),  and Sobolev's inequality, we have
\be\label{I_1.2}
I_{12}\les \|P(\rho)  u\|_{L^2}^2\les\|\rho u\|_{L^2}^2\les \|\sqrt{\rho}u\|_{L^2}^2\le C,
\ee
\beq\label{I_1.3}\begin{split}
I_{13}\les&\|\sqrt{\rho}\|_{L^3}^2\|\mathcal L^{-1}\nabla((P(\rho)-P'(\rho)\rho)\mbox{div }u)\|_{L^6}^2\\
\les&\|\nabla \mathcal L^{-1}\nabla((P(\rho)-P'(\rho)\rho)\mbox{div }u)\|_{L^2}^2
\les\|(P(\rho)-P'(\rho)\rho)\nabla u\|_{L^2}^2\\
\leq& CB_P(\|\rho\|_{L^\infty})\|\rho\|_{L^\infty}\|\nabla u\|_{L^2}^2
\le C\|\nabla u\|_{L^2}^2,
\end{split}
\eeq
where we have used the Sobolve inequality when $\Omega=\mathbb R^3$, and
both Sobolve and Poincar\'e inequalities when $\Omega$ is a bounded domain.

Putting (\ref{I_1.1-3}), (\ref{I_1.2}) and (\ref{I_1.3}) into (\ref{clcblp2.17}), and choosing $\va$ sufficiently small, we obtain
\beq\label{I_1}\begin{split}
I_1\le&
\frac{1}{2}\|\sqrt{\rho}w_t\|_{L^2}^2+C(\|\rho^{\frac{1}{5}}u\|_{L^5}^5  \|\nabla u\|_{L^{2}}^2
+\|\nabla w\|_{L^2}^2+\|\nabla u\|_{L^2}^2+\|\nabla d\|_{L^{\infty}}^2\|\de d\|_{L^2}^2+1)
\\
\le&\frac{1}{2}\|\sqrt{\rho}w_t\|_{L^2}^2+C(\|\rho^{\frac{1}{5}}u\|_{L^5}^5 \|\nabla u\|_{L^{2}}^2
+\|\nabla v\|_{L^2}^2+\|\nabla u\|_{L^2}^2+\|\nabla d\|_{L^{\infty}}^2\|\de d\|_{L^2}^2+1)
\\ 
\le&\frac{1}{2}\|\sqrt{\rho}w_t\|_{L^2}^2+C(\|\rho^{\frac{1}{5}}u\|_{L^5}^5 \|\nabla u\|_{L^{2}}^2
+\|\nabla u\|_{L^2}^2+\|\nabla d\|_{L^{\infty}}^2\|\de d\|_{L^2}^2+1),
\end{split}
\eeq
where we have used (\ref{w^{1,q} of v}) with $q=2$. 
For $I_3$, using Cauchy's inequality, we have
\beq\label{I_3}\begin{split}
I_3\le& \frac{1}{2}\int|\nabla d_t|^2\,dx+C\int |\nabla d|^2|\nabla w|^2\,dx\\ \le &\frac{1}{2}\int|\nabla d_t|^2\,dx+C\|\nabla d\|_{L^\infty}^2\int|\nabla w|^2\,dx.
\end{split}
\eeq
Substituting (\ref{I_1}) and (\ref{I_3}) into (\ref{clcblp2.10}), we obtain
\beq\label{clcblp2.11-1}
\begin{split}
&\frac{d}{dt}\int\left(\mu|\nabla w|^2+(\mu+\lambda)|\mbox{div }w|^2\right)dx+\frac{1}{2}\int\rho|w_t|^2dx\\
\le& 2\frac{d}{dt}\int (\nabla d\otimes\nabla d-\frac12{|\nabla d|^2}\mathbb I_3):\nabla w\,dx
+\frac{1}{2}\|\nabla d_t\|_{L^2}^2\\
+&C\Big(\|\nabla d\|_{L^{\infty}}^2(\|\nabla w\|_{L^2}^2+\|\de d\|_{L^2}^2)
+\|\rho^{\frac{1}{5}}u\|_{L^5}^5 \|\nabla u\|_{L^{2}}^2
+\|\nabla u\|_{L^2}^2+1\Big).
\end{split}
\eeq

{\noindent\it Step 2. Estimates of $\displaystyle\int\rho|u|^5\,dx$.}
Multiplying (\ref{clc-2}) by $5|u|^{3}u$, integrating over $\om$, and using integration by parts and Cauchy's inequality, we have
\bex
\begin{split}
&\frac{d}{dt}\int\rho|u|^5dx+\int 5|u|^{3}\left(\mu|\nabla u|^2+(\mu+\lambda)|\mbox{div }u|^2+3\mu|\nabla |u||^2\right)\,dx\\
=&\int 5 P(\rho) \mbox{div}(|u|^{3}u)\,dx+\int 5\left(\nabla d\otimes\nabla d-\frac12{|\nabla d|^2}\mathbb{I}_3\right)\mbox{div}(|u|^{3}u)
-\int 15(\mu+\lambda)(\mbox{div }u)|u|^{2}u\cdot\nabla|u|\\
\le&C(\int \rho|u|^{3}|\nabla u|+\int|\nabla d|^2|u|^{3}|\nabla u|)
+\int 5(\mu+\lambda)|u|^{3}|\mbox{div }u|^2
+\int\frac{45}{4}(\mu+\lambda)|u|^{3}\big|\nabla|u|\big|^2.
\end{split}
\eex
By Kato's inequality $|\nabla u|^2\ge |\nabla|u||^2$, we have
$$\begin{cases} 
(15\mu-\frac{45(\mu+\lambda)}{4})\int |u|^3|\nabla|u||^2\ge (15\mu-\frac{45(\mu+\lambda)}{4})\int |u|^3|\nabla u|^2,
& {\rm{if}}\ \mu-\frac{3(\mu+\lambda)}{4}\le 0\\
(15\mu-\frac{45(\mu+\lambda)}{4})\int |u|^3|\nabla|u||^2\ge 0, & {\rm{if}}\ \mu-\frac{3(\mu+\lambda)}{4}>0.
\end{cases}$$
Hence we obtain
\beq\label{rho u^5-1}
\begin{split}
&\frac{d}{dt}\int\rho|u|^5\,dx
+5\min\Big\{\mu, \ \left(4\mu-\frac{9(\mu+\lambda)}{4}\right)\Big\}\int|u|^{3}|\nabla u|^2\,dx\\
\le&C(\int\rho |u|^{3}|\nabla u|\,dx+\int|\nabla d|^2|u|^{3}|\nabla u|\,dx).
\end{split}
\eeq
Since $\lambda<\frac{7\mu}{9}$, we have
\beq\label{rho u^5-3}c_0:=5\min\Big\{\mu,\ \left(4\mu-\frac{9(\mu+\lambda)}{4}\right)\Big\}>0.\eeq
Thus by Cauchy's inequality, we have
\bex
\begin{split}
&\frac{d}{dt}\int\rho|u|^5\,dx+c_0\int|u|^{3}|\nabla u|^2\,dx\\ \le&
C(\int\rho |u|^{3}|\nabla u|\,dx+\int |\nabla d|^2|u|^{3}|\nabla u|\,dx)
\\
\le&\frac{c_0}{2}\int|u|^{3}|\nabla u|^2\,dx+C\left[\int\rho^{2}|u|^{3}\,dx+\int|\nabla d|^4|u|^{3}\,dx\right].
\end{split}
\eex
Hence by H\"older's inequality, Sobolev's inequality, the conservation of mass, (\ref{2.1}) and Young's inequality, we have
\bex
\begin{split}
&\frac{d}{dt}\int\rho|u|^5\,dx+\frac{c_0}{2}\int|u|^{3}|\nabla u|^2\,dx
\lesssim \int\rho^{2}|u|^{3}\,dx+\int|\nabla d|^4|u|^{3}\,dx\\
\lesssim&\left(\int\rho|u|^{5}\,dx\right)^{\frac{3}{5}}\left(\int\rho^{\frac72}dx\right)^{\frac{2}{5}}
+\|\nabla(|u|^{\frac{5}{2}})\|_{L^2}^{\frac{6}{5}}\left(\int|\nabla d|^5\,dx\right)^{\frac{4}{5}}\\
\le&\frac{c_0}{4}\int|u|^{3}|\nabla u|^2\,dx+C\Big[1+\int\rho|u|^5\,dx+(\int|\nabla d|^5\,dx)^2
\Big].
\end{split}
\eex
Thus by (\ref{Energy-identity}) we have
\beq\label{clcblp2.9}
\begin{split}
\frac{d}{dt}\int\rho|u|^5dx+\frac{c_0}{4}\int|u|^{3}|\nabla u|^2dx
\lesssim&\int\rho|u|^5dx
+\left(\int|\nabla d|^5dx\right)^2+1\\
\lesssim&\int\rho|u|^5dx+\|\nabla d\|_{L^{\infty}}^3\|\nabla d\|_{L^5}^5+1.
\end{split}
\eeq

{\noindent\it Step 3. Estimates of $\displaystyle\int|\nabla d|^5\,dx.$}
Differentiating (\ref{clc-3}) with respect to $x$, we obtain
\beq\label{clcblp2.12}
\nabla d_t-\nabla\de d+\nabla(u\cdot\nabla d)=\nabla(|\nabla d|^2d).
\eeq
Multiplying (\ref{clcblp2.12}) by $5|\nabla d|^3\nabla d$ and integrating by parts over $\om$, we have
\bex
\begin{split}
&\frac{d}{dt}\int|\nabla d|^5dx+5\int|\de d|^2|\nabla d|^3dx\\
=& 5\int\left[\nabla(|\nabla d|^2d)-\nabla(u\cdot\nabla d)\right]\cdot|\nabla d|^3\nabla d\,dx-5\int\de d\cdot\nabla(|\nabla d|^3)\cdot\nabla d \,dx\\
\le& \int\left(|\nabla d|^5|\nabla^2 d|+|\nabla d|^7+|\nabla u||\nabla d|^5+|\nabla d|^3|\nabla^2 d|^2\right)\,dx.
\end{split}
\eex
This, combined with Cauchy's inequality and the fact
\be\label{nabla d and Delta d}|\nabla d|^2=-d\cdot\de d\ (\mathrm{since}\ |d|=1),\ee
gives
\beq\label{clcblp2.13}
\begin{split}
&\frac{d}{dt}\int|\nabla d|^5\,dx+5\int|\de d|^2|\nabla d|^3\,dx\\
\les&\int\left(|\nabla d|^3|\nabla^2 d|^2+|\nabla d|^7+|\nabla u||\nabla d|^3|\nabla^2 d|\right)\,dx\\
\les&\|\nabla d\|_{L^{\infty}}^3\|\nabla^2 d\|_{L^{2}}^2+\|\nabla d\|_{L^{\infty}}^2\|\nabla d\|_{L^{5}}^5
+\|\nabla d\|_{L^\infty}^3\|\nabla u\|_{L^2}\|\nabla^2 d\|_{L^2}\\
\les& \|\nabla d\|_{L^{\infty}}^3(\|\nabla^2 d\|_{L^{2}}^2+\|\nabla u\|_{L^{2}}^2)+\|\nabla d\|_{L^{\infty}}^2\|\nabla d\|_{L^{5}}^5.
\end{split}
\eeq
By (\ref{w^{1,q} of v}) and (\ref{clcblp2.13}), we have
\beq\label{nabla d L^5}
\begin{split}
&\frac{d}{dt}\int|\nabla d|^5dx+5\int|\de d|^2|\nabla d|^3dx\\
\les& \|\nabla d\|_{L^{\infty}}^3(\|\nabla^2 d\|_{L^{2}}^2+\|\nabla w\|_{L^{2}}^2)+\|\nabla d\|_{L^{\infty}}^3+\|\nabla d\|_{L^{\infty}}^2\|\nabla d\|_{L^{5}}^5.
\end{split}
\eeq \\
{\noindent\it Step 4. Estimates of $\displaystyle \int|\nabla^2d|^2\,dx$.}
Multiplying (\ref{clcblp2.12}) by $\nabla d_t$, integrating by parts over $\om$, and using Cauchy's inequality, we have
\beq\label{clcblp2.14}
\begin{split}
&\frac{1}{2}\frac{d}{dt}\int |\de d|^2dx+\int |\nabla d_t|^2dx\\=&\int \left(\nabla(|\nabla d|^2d)-\nabla(u\cdot\nabla d)\right)\cdot\nabla d_tdx\\
\le&\va\|\nabla d_t\|_{L^2}^2+C\int \left(|\nabla d|^6+|\nabla d|^2|\nabla^2 d|^2+|\nabla u|^2|\nabla d|^2+|u|^2|\nabla^2 d|^2\right)dx\\
\le&\va\|\nabla d_t\|_{L^2}^2+C\Big[\|\nabla d\|_{L^{\infty}}^2\left(\|\nabla^2 d\|_{L^2}^2+\|\nabla u\|_{L^2}^2\right)+\int| u|^2|\nabla^2 d|^2\,dx\Big],
\end{split}
\eeq
where we have used (\ref{nabla d and Delta d}) to estimate
\be\label{nabla d^6}
\int |\nabla d|^6\,dx\les \|\nabla d\|_{L^{\infty}}^2\int |\nabla^2d|^2\,dx.
\ee
For the last term on the right hand side of (\ref{clcblp2.14}), using Nirenberg's interpolation inequality and Cauchy's inequality, we have
\beq\label{u^2nabla d^2-1}
\begin{split}
\int|u|^2|\nabla^2 d|^2dx\les&\left\||u|^{\frac{5}{2}}\right\|_{L^6}^{\frac{12}{15}}\|\nabla^2 d\|_{L^{\frac{30}{13}}}^2
\le\va\|\nabla |u|^{\frac{5}{2}}\|_{L^2}^2+C\|\nabla^2 d\|_{L^{\frac{30}{13}}}^{\frac{10}{3}}\\
\le&\va\|\nabla |u|^{\frac{5}{2}}\|_{L^2}^2+C\|\nabla d\|_{L^6}^{\frac{8}{3}}
\|\nabla d\|_{H^2}^{\frac23}\\
\le&\va\|\nabla |u|^{\frac{5}{2}}\|_{L^2}^2+\va\|\nabla^3 d\|_{L^2}^{2}+C(\|\nabla d\|_{L^6}^{4}+\|\nabla^2 d\|_{L^2}^{2}+1)\\
\le&5\va\int|u|^3|\nabla u|^2\,dx+\va\|\nabla^3 d\|_{L^2}^{2}+C(\|\nabla^2 d\|_{L^{2}}^{4}+1).
\end{split}\eeq
By (\ref{clc-3}), $H^3$-estimate for elliptic equations, and (\ref{nabla d^6}), we have
\beq\label{nabla^3d}\begin{split}
\|\nabla^3 d\|_{L^2}^{2}\les&\|\nabla d_t\|_{L^2}^{2}+\|\nabla(u\cdot\nabla d)\|_{L^2}^{2}+\|\nabla( |\nabla d|^2d)\|_{L^2}^{2}\\
\les&\|\nabla d_t\|_{L^2}^{2}+\|\nabla d\|_{L^{\infty}}^{2}\left(\|\nabla u\|_{L^2}^{2}+\|\nabla^2 d\|_{L^2}^{2}\right)+\|\nabla d\|_{L^6}^{6}+\int|u|^2|\nabla^2 d|^2dx\\
\le&C\Big[\|\nabla d_t\|_{L^2}^{2}+\|\nabla d\|_{L^{\infty}}^{2}\left(\|\nabla u\|_{L^2}^{2}+\|\nabla^2 d\|_{L^2}^{2}\right)+\int|u|^2|\nabla^2 d|^2\,dx\Big].
\end{split}\eeq
Substituting (\ref{nabla^3d}) into (\ref{u^2nabla d^2-1}), and choosing $\va$ sufficiently small, we have
\beq\label{u^21nabla d^2}\begin{split}
\int|u|^2|\nabla^2 d|^2dx
\le& C\Big[\int|u|^3|\nabla u|^2\,dx+\|\nabla^2 d\|_{L^{2}}^{4}+\va \|\nabla d_t\|_{L^2}^{2}\\&+
\|\nabla d\|_{L^{\infty}}^{2}\left(\|\nabla u\|_{L^2}^{2}+\|\nabla^2 d\|_{L^2}^{2}\right)+1\Big].
\end{split}\eeq
Substituting (\ref{u^21nabla d^2}) into (\ref{clcblp2.14}), using (\ref{w^{1,q} of v}), and choosing $\va$ sufficiently small, we obtain
\beq\label{clcblp2.15}
\begin{split}
& \frac{d}{dt}\int|\de d|^2dx+\int |\nabla d_t|^2dx\\
\le& C\Big[\|\nabla d\|_{L^{\infty}}^2\left(\|\nabla^2 d\|_{L^2}^2+\|\nabla u\|_{L^2}^2\right)+\int|u|^3|\nabla u|^2\,dx+\|\nabla^2 d\|_{L^{2}}^{4}+1\Big]\\
\le& C\Big[\|\nabla d\|_{L^{\infty}}^2\left(\|\nabla^2 d\|_{L^2}^2+\|\nabla w\|_{L^2}^2\right)+\|\nabla d\|_{L^{\infty}}^2+\int|u|^3|\nabla u|^2\,dx+\|\nabla^2 d\|_{L^{2}}^{4}+1\Big].
\end{split}
\eeq 

{\noindent\it Step 5. Completion of proof of Lemma \ref{clclemma2.2}.}
Adding (\ref{clcblp2.9}), (\ref{clcblp2.11-1}), (\ref{nabla d L^5}) and (\ref{clcblp2.15}) together, and choosing $\va$ sufficiently small, we obtain
\bex
\begin{split}
&\frac{d}{dt}\int\left(\rho|u|^5+\mu|\nabla w|^2+(\mu+\lambda)|\mbox{div }w|^2+|\nabla d|^5+|\de d|^2\right)dx +\frac{1}{2}\int \rho|w_t|^2dx+\frac{1}{2}\int |\nabla d_t|^2\,dx
\\
\le&2\frac{d}{dt}\int (\nabla d\otimes\nabla d-\frac12{|\nabla d|^2}\mathbb{I}_3):\nabla w\,dx
+C\Big[\left(
\|\nabla u\|_{L^{2}}^2+1\right)\int \rho|u|^5dx\|\nabla u\|_{L^2}^2+\|\nabla d\|_{L^{\infty}}^3\\
&+\|\nabla d\|_{L^{\infty}}^3\left(\|\nabla d\|_{L^5}^5+\|\nabla^2 d\|_{L^{2}}^2+\|\nabla w\|_{L^{2}}^2\right)+\|\nabla d\|_{L^{\infty}}^2\left(\|\nabla d\|_{L^{5}}^5+\|\nabla^2 d\|_{L^2}^2+\|\nabla w\|_{L^2}^2\right)
\\&+\|\nabla d\|_{L^{\infty}}^2+\|\nabla^2 d\|_{L^{2}}^{4}+1\Big].
\end{split}
\eex
This, combined with Cauchy's inequality, implies
\bex
\begin{split}
&\frac{d}{dt}\int\left(\rho|u|^5+\mu|\nabla w|^2+(\mu+\lambda)|\mbox{div }w|^2+|\nabla d|^5+|\de d|^2\right)dx +\frac{1}{2}\left(\int \rho|w_t|^2\,dx+\int |\nabla d_t|^2\,dx\right)
\\
\le&2\frac{d}{dt}\int (\nabla d\otimes\nabla d-\frac12{|\nabla d|^2}\mathbb{I}_3):\nabla w\,dx
+C\Big[\|\nabla u\|_{L^2}^2+\|\nabla d\|_{L^{\infty}}^3\\
&+\left(\|\nabla u\|_{L^{2}}^2+\|\nabla^2 d\|_{L^{2}}^2+\|\nabla d\|_{L^{\infty}}^3+1\right)\left(\|\rho^\frac{1}{5}u\|_{L^5}^5+\|\nabla d\|_{L^{5}}^5+\|\nabla^2 d\|_{L^2}^2+\|\nabla w\|_{L^2}^2\right)+1\Big].
\end{split}
\eex
Integrating over $(0,t)$, and using (\ref{2.1}), (\ref{Energy-identity}), we have
\beq\label{total-1}
\begin{split}
&\int\left(\rho|u|^5+|\nabla w|^2+|\nabla d|^5+|\nabla^2 d|^2\right)dx +\int_0^t\int \left(\rho|w_t|^2+|\nabla d_t|^2\right)\,dx\,ds
\\
\le&C\Big[\int |\nabla d|^2|\nabla w|\,dx
+\int_0^tK(s)\left(\|\rho^\frac{1}{5}u\|_{L^5}^5+\|\nabla w\|_{L^2}^2+\|\nabla d\|_{L^{5}}^5+\|\nabla^2 d\|_{L^2}^2\right)\,ds+1\Big],
\end{split}
\eeq
where $$K(s)=\|\nabla u(s)\|_{L^{2}}^2+\|\nabla^2 d(s)\|_{L^{2}}^2+\|\nabla d(s)\|_{L^{\infty}}^3+1.$$
By (\ref{total-1}) and Young's inequality, we have
\bex
\begin{split}
&\int\left(\rho|u|^5+|\nabla w|^2+|\nabla d|^5+|\nabla^2 d|^2\right)dx +\int_0^t\int \left(\rho|w_t|^2+|\nabla d_t|^2\right)\,dx\,ds
\\
\le&\frac{1}{2}\int|\nabla w|^2\,dx+C\Big[\int|\nabla d|^4\,dx
+\int_0^tK(s)\left(\|\rho^\frac{1}{5}u\|_{L^5}^5+\|\nabla w\|_{L^2}^2+\|\nabla d\|_{L^{5}}^5+\|\nabla^2 d\|_{L^2}^2\right)\,ds+1\Big]
\\
\le&\frac{1}{2}(\int|\nabla w|^2\,dx+\int|\nabla d|^5\,dx)
+C\Big[\int_0^tK(s)\left(\|\rho^\frac{1}{5}u\|_{L^5}^5+\|\nabla w\|_{L^2}^2+\|\nabla d\|_{L^{5}}^5+\|\nabla^2 d\|_{L^2}^2\right)\,ds+1\Big].
\end{split}
\eex
Thus we obtain
\beq\label{total-2}
\begin{split}
&\int\left(\rho|u|^5+|\nabla w|^2+|\nabla d|^5+|\nabla^2 d|^2\right)dx +\int_0^t\int \left(\rho|w_t|^2+|\nabla d_t|^2\right)\,dx\,ds
\\
\le&C\Big[1+\int_0^tK(s)\left(\|\rho^\frac{1}{5}u\|_{L^5}^5+\|\nabla w\|_{L^2}^2+\|\nabla d\|_{L^{5}}^5+\|\nabla^2 d\|_{L^2}^2\right)\,ds\Big].
\end{split}
\eeq
By (\ref{2.1}), (\ref{Energy-identity}) and (\ref{clcblp2.2}), we know
\beq\label{total-3}
\int_0^tK(s)ds\le C.
\eeq
By (\ref{total-2}), (\ref{total-3}) and Gronwall's inequality, we obtain that for any $0\le t<T_*$,
\bex
\int\left(\rho|u|^5+|\nabla w|^2+|\nabla d|^5+|\nabla^2 d|^2\right)dx +\int_0^t\int \left(\rho|w_t|^2+|\nabla d_t|^2\right)\,dx\,ds\le C.
\eex 
This completes the proof of Lemma \ref{clclemma2.2}.
\endpf

\begin{corollary}\label{clccoro2.3}{\it  Under the same assumptions of Lemma \ref{clclemma2.2}, we have
that for any $2\le q\le 6$,
\beq\label{clcblp3.16}
\sup\ls_{0\leq t<T_*}\left(\|u\|_{L^{6}}+\|\nabla u\|_{L^{2}}+\|\nabla d\|_{L^{q}}+\|d_t\|_{L^{2}}\right)+\|\nabla u\|_{L^{2}(0,T;L^6)}\leq C.
\eeq
}
\end{corollary}

\pf Combining (\ref{w^{1,q} of v}) with (\ref{clcblp2.8}), we get
\be\label{nabla u L^2}\|\nabla u(t)\|_{L^2}\les \|\nabla w(t)\|_{L^2}+\|\nabla v(t)\|_{L^2}\leq C.\ee
The upper bound of $\sup\ls_{0\leq t<T_*}\|u\|_{L^{6}}$ 
follows from  (\ref{nabla u L^2}) and Sobolev's inequality. The bound of $\sup\ls_{0\leq t<T_*}\|\nabla d\|_{L^{q}}$ follows from  (\ref{clcblp2.8}) and interpolation inequality.
For the last term of (\ref{clcblp3.16}), by Sobolev's inequality, (\ref{w^{1,q} of v}) and (\ref{clcblp2.8}), we have
\bex\begin{split}
\|\nabla u\|_{L^{2}(0,T;L^6)}\les&\|\nabla w\|_{L^{2}(0,T;L^6)}+\|\nabla v\|_{L^{2}(0,T;L^6)}\\
\les&\|\nabla^2 w\|_{L^{2}(0,T;L^2)}+\|\nabla w\|_{L^{2}(0,T;L^2)}+1\leq C.
\end{split}
\eex
By equation (\ref{clc-3}), (\ref{clcblp2.8}) and H\"older's inequality, we have
\bex\begin{split}
\sup\ls_{0\leq t<T_*}\|d_t\|_{L^{2}}
\les&\sup\ls_{0\leq t<T_*}\left(\|\de d\|_{L^{2}}+\|\nabla d\|_{L^{4}}^2+\|u\cdot\nabla d\|_{L^{2}}\right)\\
\les&\sup\ls_{0\leq t<T_*}\left(\|u\|_{L^{6}}\|\nabla d\|_{L^{3}}\right)+1\leq C.
\end{split}
\eex
This completes the proof. \endpf

\blm\label{clclemma3.1}{\it Under the same assumptions of Lemma \ref{clclemma2.2}, 
 $(\rho,u,d)$ satisfies that for any $0\leq t<T_*$,
\beq\label{clcblp3.1}
\begin{split}
\int_{\om}(\rho|\dot{u}(t)|^2+|\nabla d_t|^2)(t)\,dx+\int_0^t\int_{\om}(|\nabla \dot{u}|^2+|d_{tt}|^2)\,dxds\leq C,
\end{split}
\eeq
where $\dot{f}$ is the material derivative:
$$\dot{f}:=f_t+u\cdot\nabla f.$$}
\elm

\pf {\noindent\it Step 1. Estimates of $\displaystyle\int\rho|\dot{u}(t)|^2\,dx$.}
By the definition of material derivative, we can write (\ref{clc-2}) as follows,
\beq\label{clc-2-1}\rho\dot{u}+\nabla (P(\rho))=\mathcal Lu-\nabla d\cdot\de d.\eeq
Differentiating (\ref{clc-2-1}) with respect to $t$ and using (\ref{clc-1}), we have
\beq\label{clcblp3.2}
\begin{split}
&\rho\u_t+\rho u\cdot\nabla\u+\nabla (P(\rho)_t)+(\nabla d\cdot\de d)_t\\
=&\mathcal L\u-\mathcal L(u\cdot\nabla u)+\di\Big[\mathcal Lu\otimes u-\nabla (P(\rho))\otimes u-(\nabla d\cdot\de d)\otimes u\Big].
\end{split}
\eeq
Multiplying (\ref{clcblp3.2}) by $\u$, integrating by parts over $\om$ and using the fact $\u=0$ on $\partial\om$, we obtain
\beq\label{clcblp3.3}
\begin{split}
&\frac{1}{2}\frac{d}{dt}\int\rho|\u|^2\,dx+\int\left(\mu|\nabla\u|^2+(\mu+\lambda)|\di\u|^2\right)\,dx\\
=&\int\left( (P(\rho))_t\di\u + u\otimes\nabla (P(\rho)):\nabla\u\right) \,dx+\mu\int\left(\di(\de u\otimes u)-\de(u\cdot\nabla u)\right)\cdot\u \,dx\\
+&(\mu+\lambda)\int\Big(\di(\nabla\di u\otimes u)-\nabla\di(u\cdot\nabla u)\Big)\cdot\u \,dx
+\int(u\otimes (\de d\cdot\nabla d)):\nabla\u \,dx\\
+&\int(\nabla d_t\otimes\nabla d+\nabla d\otimes\nabla d_t-\nabla d\cdot\nabla d_t\mathbb{I}_3):\nabla\u \,dx
=\sum\ls_{i=1}^5J_i.
\end{split}
\eeq
By equation (\ref{clc-1}) and (\ref{2.1}), we have
\bex
\begin{split}
&J_1=\int\Big(-\di(P(\rho)u)\di\u -(P'(\rho)\rho-P(\rho))\di u\di\u + u\otimes\nabla (P(\rho)):\nabla\u \Big)\,dx\\
=&\int\Big(P(\rho)u\cdot\nabla\di\u + (P(\rho)-P'(\rho)\rho)\di u\di\u +P(\rho)(\nabla u)^t:\nabla\u -P(\rho) u\cdot\nabla\di\u\Big) \,dx\\
=&\int\Big((P(\rho)-P'(\rho)\rho)\di u\di\u \,dx+ P(\rho)(\nabla u)^t:\nabla\u \Big)\,dx
\les\|\nabla u\|_{L^2}\|\nabla \u\|_{L^2}.
\end{split}
\eex
By the product rule, we can see 
$$\di(\de u\otimes u)-\de(u\cdot\nabla u)
=\nabla_k(\di u\nabla_ku)-\nabla_k(\nabla_ku^j\nabla_ju)-\nabla_j(\nabla_ku^j\nabla_ku),$$
so that by integration by parts, we have
\bex
\begin{split}
J_2
=\mu\int\left(\nabla_k(\di u\nabla_ku)-\nabla_k(\nabla_ku^j\nabla_ju)-\nabla_j(\nabla_ku^j\nabla_ku)\right)\cdot\u \,dx
\les\|\nabla \u\|_{L^2}\|\nabla u\|_{L^4}^2.
\end{split}
\eex
Similarly, since
$$
\di(\nabla\di u\otimes u)-\nabla\di(u\cdot\nabla u)=\nabla_k(\nabla_ju^j\nabla_iu^i)-\nabla_k(\nabla_ju^i\nabla_iu^j)-\nabla_i(\nabla_ku^i\nabla_ju^j),
$$
we have
\bex
\begin{split}
J_3=
(\mu+\lambda)\int\left(\nabla_k(\nabla_ju^j\nabla_iu^i)-\nabla_k(\nabla_ju^i\nabla_iu^j)-\nabla_i(\nabla_ku^i\nabla_ju^j)\right)\u^k \,dx
\les\|\nabla \u\|_{L^2}\|\nabla u\|_{L^4}^2.
\end{split}
\eex
By H\"older's inequality, and Corollary \ref{clccoro2.3}, we have
\bex
\begin{split}
J_4\les\|\nabla\u\|_{L^2}\|\de d\|_{L^6}\|\nabla d\|_{L^6}\|u\|_{L^6}
\les\|\nabla\u\|_{L^2}\|\de d\|_{L^6},
\end{split}
\eex
\bex
\begin{split}
J_5\les\int|\nabla\u| |\nabla d_t| |\nabla d|\,dx
\les\|\nabla\u\|_{L^2}\|\nabla d_t\|_{L^2}\|\nabla
d\|_{L^{\infty}}.
\end{split}
\eex
Putting all these estimates into (\ref{clcblp3.3}), using Young's inequality and Sobolev's inequality,and Lemma 4.2 and
Corollary 4.2,  we have
\bex
\begin{split}
&\frac{1}{2}\frac{d}{dt}\int\rho|\u|^2\,dx+\int\left(\mu|\nabla\u|^2+(\mu+\lambda)|\di\u|^2\right)\,dx\\
\les&\|\nabla u\|_{L^2}\|\nabla \u\|_{L^2}+\|\nabla u\|_{L^4}^2\|\nabla \u\|_{L^2}+\|\nabla\u\|_{L^2}\|\de d\|_{L^6}+\|\nabla\u\|_{L^2}\|\nabla d_t\|_{L^2}\|\nabla d\|_{L^{\infty}}\\
\leq&\frac{\mu}{2}\|\nabla \u\|_{L^2}^2+C\left(\|\nabla u\|_{L^2}^2+\|\nabla u\|_{L^4}^4+\|\de d\|_{H^1}^2+\|\nabla d_t\|_{L^2}^2\|\nabla d\|_{L^{\infty}}^2\right)\\
\leq&\frac{\mu}{2}\|\nabla \u\|_{L^2}^2+C\left(\|\nabla u\|_{L^4}^4+\|\nabla^3 d\|_{L^2}^2+\|\nabla d_t\|_{L^2}^2\|\nabla d\|_{L^{\infty}}^2+1\right)
\end{split}
\eex
Thus we obtain
\beq\label{clcblp3.4}
\begin{split}
\frac{d}{dt}\int\rho|\u|^2\,dx+\mu\int|\nabla\u|^2\,dx
\les\|\nabla u\|_{L^4}^4+\|\nabla^3 d\|_{L^2}^2+\|\nabla d_t\|_{L^2}^2\|\nabla d\|_{L^{\infty}}^2+1.
\end{split}
\eeq
By $H^3$-estimate of elliptic equations, Lemma \ref{clclemma2.2}, Corollary \ref{clccoro2.3}, and Nirenberg's interpolation inequality, we have
\bex\begin{split}
\|\nabla^3 d\|_{L^2}\les&\|\nabla d_t\|_{L^2}+\|u\cdot\nabla d\|_{L^2}+\||\nabla u| |\nabla d|\|_{L^2}+\||\nabla d| |\nabla^2 d|\|_{L^2}+\||\nabla d|^3\|_{L^2}\\
\les&\|\nabla d_t\|_{L^2}+\|u\|_{L^6}\|\nabla d\|_{L^3}+\|\nabla u\|_{L^3}\|\nabla d\|_{L^6}+\|\nabla d\|_{L^6}\|\nabla^2 d\|_{L^3}+1\\
\les&\|\nabla d_t\|_{L^2}+\|\nabla u\|_{L^2}^{\frac{1}{2}}\|\nabla u\|_{H^1}^\frac{1}{2}+\|\nabla^2 d\|_{L^2}^{\frac{1}{2}}\|\nabla^2 d\|_{H^1}^\frac{1}{2}+1\\
\leq&\frac{1}{2}\|\nabla^3d\|_{L^2}+C\left(\|\nabla d_t\|_{L^2}+\|\nabla^2 u\|_{L^2}+1\right).
\end{split}
\eex
Thus we obtain
\beq\label{clcblp3.10}
\|\nabla^3 d\|_{L^2}\les\|\nabla d_t\|_{L^2}+\|\nabla^2 u\|_{L^2}+1.
\eeq
By the definition of $w$, we have
\beq\label{clcblp3.8}
\mathcal Lw=\rho\u+\de d\cdot\nabla d.
\eeq
By $H^2$-estimate of the equation (\ref{clcblp3.8}), (\ref{2.1}), Corollary \ref{clccoro2.3}, Nirenberg's interpolation inequality, and (\ref{clcblp3.10}),  we obtain
\beq\label{clcblp3.11}
\begin{split}
\|\nabla^2 w\|_{L^2}^2\les& \|\rho\u\|_{L^2}^2+\|\de d\cdot\nabla d\|_{L^2}^2\les \|\rho^{\frac{1}{2}}\u\|_{L^2}^2+\|\de d\|_{L^3}^2\|\nabla d\|_{L^6}^2\\
\les& \|\rho^{\frac{1}{2}}\u\|_{L^2}^2+\|\de d\|_{L^2}\|\de d\|_{H^1}\les \|\rho^{\frac{1}{2}}\u\|_{L^2}^2+\|\nabla\de d\|_{L^2}+1\\
\les& \|\rho^{\frac{1}{2}}\u\|_{L^2}^2+\|\nabla d_t\|_{L^2}+\|\nabla^2 u\|_{L^2}+1.
\end{split}
\eeq
By interpolation inequality, Corollary \ref{clccoro2.3}, (\ref{w^{1,q} of v}) (for $q=6$), (\ref{clcblp3.11}), and Cauchy's inequality, we obtain
\beq\label{clcblp3.12}
\begin{split}
\|\nabla u\|_{L^4}^4\les&\|\nabla u\|_{L^2}\|\nabla u\|_{L^6}^3\les\|\nabla u\|_{L^6}\|\nabla u\|_{L^6}^2\\
\les& \|\nabla u\|_{L^6}\left(\|\nabla w\|_{L^6}^2+\|\nabla v\|_{L^6}^2\right)
\les \|\nabla u\|_{L^6}\left(\|\nabla^2 w\|_{L^2}^2+1\right)\\
\les& \|\nabla u\|_{L^6}\left(\|\rho^{\frac{1}{2}}\u\|_{L^2}^2+\|\nabla d_t\|_{L^2}+\|\nabla^2 u\|_{L^2}+1\right)\\
\les& \|\nabla u\|_{L^6}\|\rho^{\frac{1}{2}}\u\|_{L^2}^2+\|\nabla u\|_{L^6}^2+\|\nabla d_t\|_{L^2}^2+\|\nabla^2 u\|_{L^2}^2+1\\
\les& \|\nabla u\|_{L^6}\|\rho^{\frac{1}{2}}\u\|_{L^2}^2+\|\nabla d_t\|_{L^2}^2+\|\nabla^2 u\|_{L^2}^2+1.
\end{split}
\eeq
Putting (\ref{clcblp3.12}) and (\ref{clcblp3.10}) into (\ref{clcblp3.4}), we have
\beq\label{clcblp3.5}
\begin{split}
\frac{d}{dt}\int\rho|\u|^2\,dx+\mu\int|\nabla\u|^2\,dx
\les \|\nabla u\|_{L^6}\|\rho^{\frac{1}{2}}\u\|_{L^2}^2+\|\nabla d_t\|_{L^2}^2(\|\nabla d\|_{L^{\infty}}^2+1)+\|\nabla^2 u\|_{L^2}^2+1.
\end{split}
\eeq
{\noindent\it Step 2. Estimates of $\displaystyle\int |\nabla d_t|^2\,dx$.}
Differentiating (\ref{clc-3}) with respect to $t$, we have
\beq\label{clcblp3.13}
d_{tt}-\de d_t=\pa_t\left(|\nabla d|^2d-u\cdot\nabla d\right).
\eeq
Multiplying (\ref{clcblp3.13}) by $d_{tt}$, integrating by parts over $\om$ and using $\frac{\partial d_t}{\partial \nu}\Big|_{\partial\om}=0$, we obtain
\beq\label{clcblp3.14}\begin{split}
&\frac{1}{2}\frac{d}{dt}\int|\nabla d_t|^2\,dx+\int |d_{tt}|^2\,dx=\int\pa_t\left(|\nabla d|^2d-u\cdot\nabla d\right)d_{tt}\,dx\\
\les&\int\left(|\nabla d|^2|d_t|+|\nabla d||\nabla d_t|\right)|d_{tt}|\,dx
+\int\left(|u_t||\nabla d|+|u||\nabla d_t|\right)|d_{tt}|\,dx
=K_1+K_2.
\end{split}
\eeq
By H\"older's inequality, Sobolev's inequality, Corollary \ref{clccoro2.3}, and Young's inequality, we have
\bex\begin{split}
|K_1|\les&\|d_{tt}\|_{L^2}\|d_t\|_{L^6}\|\nabla d\|_{L^6}^2+\|d_{tt}\|_{L^2}\|\nabla d_t\|_{L^2}\|\nabla d\|_{L^{\infty}}\\
\les&\|d_{tt}\|_{L^2}(\|\nabla d_t\|_{L^2}+1)+\|d_{tt}\|_{L^2}\|\nabla d_t\|_{L^2}\|\nabla d\|_{L^{\infty}}\\
\leq&\frac{1}{8}\|d_{tt}\|_{L^2}^2+C\left(\|\nabla d_t\|_{L^2}^2+\|\nabla d\|_{L^{\infty}}^2\|\nabla d_t\|_{L^2}^2+1\right).
\end{split}
\eex
By the definition of $\u$, H\"older's inequality, Sobolev's inequality, Corollary \ref{clccoro2.3}, and Young's inequality, we have
\bex\begin{split}
|K_2|\les&\int\left[(|\u|+|u||\nabla u|)|\nabla d|+|u||\nabla d_t|\right]|d_{tt}|\,dx\\
\les&\|d_{tt}\|_{L^2}\|\u\|_{L^6}\|\nabla d\|_{L^3}+\|d_{tt}\|_{L^2}\|u\|_{L^6}\|\nabla u\|_{L^6}\|\nabla d\|_{L^6}
+\|d_{tt}\|_{L^2}\|u\|_{L^6}\|\nabla d_t\|_{L^3}\\
\les&\|d_{tt}\|_{L^2}\|\nabla\u\|_{L^2}+\|d_{tt}\|_{L^2}(\|\nabla^2 u\|_{L^2}+1)+\|d_{tt}\|_{L^2}\|\nabla d_t\|_{L^3}\\
\leq&\frac{1}{8}\|d_{tt}\|_{L^2}^2+C\left(\|\nabla\u\|_{L^2}^2+\|\nabla^2 u\|_{L^2}^2+\|\nabla d\|_{L^3}^2+1\right).
\end{split}
\eex
Putting these two estimates into (\ref{clcblp3.14}), using Nirenberg's interpolation inequality, and Young's inequality, we have
\beq\label{clcblp3.6}
\begin{split}
&\frac{1}{2}\frac{d}{dt}\int|\nabla d_t|^2\,dx+\frac{3}{4}\int|d_{tt}|^2\,dx\\
\les&\|\nabla\u\|_{L^2}^2+\|\nabla^2 u\|_{L^2}^2+(1+\|\nabla d\|_{L^{\infty}}^2)\|\nabla d_t\|_{L^2}^2+\|\nabla d\|_{L^3}^2+1\\
\les&\|\nabla\u\|_{L^2}^2+\|\nabla^2 u\|_{L^2}^2+(1+\|\nabla d\|_{L^{\infty}}^2)\|\nabla d_t\|_{L^2}^2+\|\nabla d_t\|_{L^2}\|\nabla^2 d_t\|_{L^2}+1\\
\leq&\frac{1}{8}\|\nabla^2 d_t\|_{L^2}^2+C\left(\|\nabla\u\|_{L^2}^2+\|\nabla^2 u\|_{L^2}^2+(1+\|\nabla d\|_{L^{\infty}}^2)\|\nabla d_t\|_{L^2}^2+1\right).
\end{split}
\eeq
By $H^2$-estimate of the equation (\ref{clcblp3.13}) and estimates similar to $K_1$ and $K_2$, we obtain
\bex\begin{split}
\|\nabla^2 d_t\|_{L^2}\les&\| d_{tt}\|_{L^2}+\|\pa_t(u\cdot\nabla d)\|_{L^2}+\|\pa_t(|\nabla d|^2d)\|_{L^2}\\
\les&\| d_{tt}\|_{L^2}+\|\u\cdot\nabla d\|_{L^2}+\|(u\cdot\nabla u)\cdot\nabla d\|_{L^2}+\|u\|_{L^6}\|\nabla d_t\|_{L^3}\\
&+\|d_t\|_{L^6}\|\nabla d\|_{L^6}^{2}+\|\nabla d_t\|_{L^3}\|\nabla d\|_{L^6}\\
\les&\| d_{tt}\|_{L^2}+\|\u\|_{L^6}\|\nabla d\|_{L^3}+\|u\|_{L^6}\|\nabla u\|_{L^6}\|\nabla d\|_{L^6}+\|\nabla d_t\|_{L^2}^{\frac{1}{2}}\|\nabla d_t\|_{L^6}^{\frac{1}{2}}+\|\nabla d_t\|_{L^2}+1\\
\leq&\frac{1}{2}\|\nabla^2 d_t\|_{L^2}+C\left(\| d_{tt}\|_{L^2}+\|\nabla\u\|_{L^2}+\|\nabla^2 u\|_{L^2}+\|\nabla d_t\|_{L^2}+1\right).
\end{split}
\eex
Thus
\beq\label{nabla^2d_t}\begin{split}
\|\nabla^2 d_t\|_{L^2}
\les&\| d_{tt}\|_{L^2}+\|\nabla\u\|_{L^2}+\|\nabla^2 u\|_{L^2}+\|\nabla d_t\|_{L^2}+1.
\end{split}
\eeq
Substituting this inequality into (\ref{clcblp3.6}), we obtain
\beq\label{clcblp3.7}
\begin{split}
\frac{d}{dt}\int|\nabla d_t|^2\,dx+\int|d_{tt}|^2\,dx\les&\|\nabla\u\|_{L^2}^2+\|\nabla^2 u\|_{L^2}^2+(1+\|\nabla d\|_{L^{\infty}}^2)\|\nabla d_t\|_{L^2}^2+1.
\end{split}
\eeq

Combining (\ref{clcblp3.5}) and (\ref{clcblp3.7}), and applying Gronwall's inequality, 
we establish the conclusions of Lemma \ref{clclemma3.1}.
\endpf
\vspace{5mm}

By the equation (\ref{clcblp3.8}) and Lemma \ref{clclemma3.1}, we obtain the following Corollary.
\begin{corollary}\label{clccoro3.2}{\it  Under the same assumptions of Lemma \ref{clclemma2.2}, we have
that for $q\in (3,6]$,
\beq\label{H^3-d}\sup\limits_{0\leq t<T_*}\left(\|\nabla^3 d\|_{L^2}+\|\nabla d\|_{L^{\infty}}\right)+\|\nabla w\|_{L^{2}(0,T_*;L^{\infty})}+\|\nabla^2 w\|_{L^{2}(0,T_*;L^q)}\leq C.
\eeq}
\end{corollary}

\pf By $H^3$-estimate of elliptic equations, (\ref{clc-3}), Lemma \ref{clclemma3.1}, Corollary \ref{clccoro2.3}, and Nirenberg's interpolation inequality, we have
\bex\begin{split}
\|\nabla^3 d\|_{L^2}\les&\|\nabla d_t\|_{L^2}+\|u\cdot\nabla d\|_{L^2}+\||\nabla u| |\nabla d|\|_{L^2}+\||\nabla d| |\nabla^2 d|\|_{L^2}+\||\nabla d|^3\|_{L^2}\\
\les&\|u\|_{L^6}\|\nabla d\|_{L^3}+\|\nabla u\|_{L^2}\|\nabla d\|_{L^{\infty}}+\|\nabla d\|_{L^6}\|\nabla^2 d\|_{L^3}+1\\
\les&\|\nabla d\|_{L^2}^{\frac{1}{4}}\|\nabla d\|_{H^2}^\frac{3}{4}+\|\nabla^2 d\|_{L^2}^{\frac{1}{2}}\|\nabla^2 d\|_{H^1}^\frac{1}{2}+1\
\leq\frac{1}{2}\|\nabla^3d\|_{L^2}+C.
\end{split}
\eex
Hence
$$\sup\limits_{0\leq t<T_*}\|\nabla^3 d\|_{L^2}\leq C.$$
By Sobolev's inequality, this yields
$$\sup\limits_{0\leq t<T_*}\|\nabla d\|_{L^{\infty}}\leq C.$$
 For simplicity, we only consider the case $q=6$. 
By $W^{2,q}$-estimate of the equation (\ref{clcblp3.8}), (\ref{2.1}), and Sobolev's inequality, we obtain
\bex\begin{split}
\|\nabla^2 w\|_{L^{6}}\les&\|\rho\u\|_{L^{6}}+\|\de d\cdot\nabla d\|_{L^{6}}
\les\|\u\|_{L^{6}}+\|\de d\|_{H^{1}}\|\nabla d\|_{L^{\infty}}\
\les\|\nabla\u\|_{L^{2}}+1.
\end{split}
\eex
Therefore, by (\ref{clcblp3.1}), we have
\bex\begin{split}
\|\nabla^2 w\|_{L^{2}(0,T_*;L^6)}\les\int_{0}^{T_*}\left(\|\nabla\u\|_{L^{2}}^2+1\right)\,ds\leq C.
\end{split}
\eex
\endpf

Following the same argument of \cite{Sun-Wang-Zhang} Section 5, we have
\begin{lemma}\label{le:2.6}{\it  Under the same assumptions of Lemma \ref{clclemma2.2}, we have
that for $q\in (3, 6]$,
\beq\label{rho-q-est}\sup\limits_{0\leq t<T_*}\|\nabla \rho\|_{L^q\bigcap L^2}\leq C.
\eeq}
\end{lemma}

\begin{corollary}\label{clccoro2.6}{\it  Under the same assumptions of Lemma \ref{clclemma2.2}, we have
for $q\in (3,6]$,
\beq\label{H^2-u-q}
\sup\limits_{0\leq t<T_*}\|\nabla^2 u\|_{L^2}+\|u\|_{L^2(0,T_*;D^{2,q})}\leq C.
\eeq}
\end{corollary}
\pf
By Proposition 2.1 in \cite{Sun-Wang-Zhang}, (\ref{clc-2-1}), (\ref{2.1}) and Lemma \ref{le:2.6}, we 
obtain that 
for $r_1=2$ or $q$,
\beq\label{nabla^2u:r_1}
\begin{split}
\|\nabla^2 u\|_{L^{r_1}}\les& \|\rho \u\|_{L^{r_1}}+\|\nabla (P(\rho))\|_{L^{r_1}}+\|\nabla d\cdot\de d\|_{L^{r_1}}\\
\les&\|\rho \u\|_{L^{r_1}}+\|\nabla d\cdot\de d\|_{L^{r_1}}+\|\nabla\rho\|_{L^{r_1}}.
 \end{split}
 \eeq 
 When $r_1=2$, (\ref{2.1}), (\ref{nabla^2u:r_1}), Lemma \ref{clclemma2.2}, Lemma \ref{clclemma3.1} ,
and Corollary \ref{clccoro3.2} imply
 \bex
\begin{split}
\|\nabla^2 u\|_{L^2}\les\|\rho^\frac{1}{2}\u\|_{L^2}+\|\nabla d\|_{L^\infty}\|\de d\|_{L^2}+1\le C.
 \end{split}
 \eex
 When $r_1=q$, for simplicity, we only consider the case $q=6$. By (\ref{2.1}), (\ref{nabla^2u:r_1}), Lemma \ref{clclemma2.2}, Lemma \ref{clclemma3.1}, Corollary \ref{clccoro3.2}, and
 Sobolev's inequality, we have
 \bex
\begin{split}
\|\nabla^2 u\|_{L^2(0,T_*;L^6)}\les& \|\rho\|_{L^\infty(0,T_*; L^\infty)}\|\u\|_{L^2(0,T_*;L^6)}+\sup\limits_{0\le t<T_*}\|\nabla d\|_{L^\infty}\|\de d\|_{L^2(0,T_*;L^6)}+1\\
\les&\|\nabla\u\|_{L^2(0,T_*;L^2)}+\|\de d\|_{L^2(0,T_*;H^1)}+1\le C.
 \end{split}
 \eex
This completes the proof. \endpf
\begin{corollary}\label{clccoro2.7}{\it  Under the same assumptions of Lemma \ref{clclemma2.2}, we have
that for $r_1=2$ or $q$, 
\beq\label{rhout-est}\sup\limits_{0\leq t<T_*}\int_\Omega\big(\rho|u_t|^2+|\rho_t|^{r_1}\big)\,dx+\int_0^{T_*}\int_\Omega\left(|\nabla u_t|^2+|\nabla^2d_t|^2+|\nabla^4d|^2\right)\,dxds\le C.
\eeq}
\end{corollary}
\pf It follows from (\ref{2.1}), Lemma \ref{clclemma3.1}, Sobolev's inequality, (\ref{clcblp3.16}),
 and Corollary \ref{clccoro2.6} that
\bex\begin{split}
\int\rho|u_t|^2\,dx\les& \int\rho|\u|^2\,dx+\int\rho|u\cdot\nabla u|^2\,dx\\ \les&\|\rho\|_{L^\infty}\|u\|_{L^\infty}\int|\nabla u|^2\,dx+1\ \les
\|\nabla u\|_{H^1}+1\le C.
\end{split}
\eex By (\ref{clc-1}), (\ref{2.1}), Sobolev's inequality, (\ref{clcblp3.16}), Lemma \ref{le:2.6} and Corollary \ref{clccoro2.6}, we get
\bex\begin{split}
\|\rho_t\|_{L^{r_1}}\les& \|\rho\mathrm{div}u\|_{L^{r_1}}+\|u\cdot\nabla\rho\|_{L^{r_1}}\les
\|\rho\|_{L^\infty}\|\mathrm{div}u\|_{L^{r_1}}+\|u\|_{L^\infty}\|\nabla\rho\|_{L^{r_1}}\\ \les&\|\rho\|_{L^\infty}\|\mathrm{div}u\|_{H^1}+\|\nabla u\|_{H^1}\|\nabla\rho\|_{L^{r_1}}\le C.
\end{split}
\eex
 By Lemma \ref{clclemma3.1}, interpolation inequality, Sobolev's inequality, (\ref{clcblp3.16}),
 and Corollary \ref{clccoro2.6}, we have
\bex
\begin{split}
\int_0^{T_*}\int_\om|\nabla u_t|^2\,dxds\les &\int_0^{T_*}\int_\om|\nabla \u|^2\,dxds+
\int_0^{T_*}\int_\om|\nabla(u\cdot\nabla u)|^2\,dxds\\ \les&
\int_0^{T_*}\int_\om|\nabla u|^4\,dxds+\int_0^{T_*}\int_\om|u\cdot\nabla^2 u|^2\,dxds+1\\ \les &\int_0^{T_*}\|\nabla u\|_{L^2}\|\nabla u\|_{H^1}^3\,ds+\int_0^{T_*}\|u\|_{L^\infty}^2\int_\om|\nabla^2 u|^2\,dxds+1\\ \les&\int_0^{T_*}\|\nabla u\|_{H^1}^2\int_\om|\nabla^2 u|^2\,dxds+1\le C.
\end{split}
\eex
By (\ref{nabla^2d_t}), Lemma \ref{clclemma3.1}, and Corollary \ref{clccoro2.6}, we get
\beq\label{nabla^2d_t-L^2}
\int_0^{T_*}\int_\om|\nabla^2d_t|^2\,dxds\le C.
\eeq
By $H^4$-estimate of the equation (\ref{clc-3}), we have
\beq\label{nabla^4d}\begin{split}
\|\nabla^4d\|_{L^2}^2\les\|d_t\|_{H^2}^2+\|u\cdot\nabla d\|_{H^2}^2+\||\nabla d|^2d\|_{H^2}^2=\sum\limits_{i=1}^3L_i.
\end{split}
\eeq
For $L_1$, (\ref{clcblp3.16}) and Lemma \ref{clclemma3.1} imply
\beq\label{nabla^4d-1}
L_1\les \|\nabla^2 d_t\|_{L^2}^2+1.
\eeq
For $L_2$, H\"older's inequality, Sobolev's inequality, (\ref{Energy-identity}), (\ref{clcblp2.8}), (\ref{clcblp3.16}), Corollary \ref{clccoro3.2}, and Corollary \ref{clccoro2.6}, we have
\beq\label{nabla^4d-2}
\begin{split}
L_2\les &\Big\||u|(|\nabla d|+|\nabla^2d|+|\nabla^3 d|)\Big\|_{L^2}^2+\Big\||\nabla u|(|\nabla d|+|\nabla^2 d|)\Big\|_{L^2}^2
+\Big\||\nabla^2 u||\nabla d|\Big\|_{L^2}^2\\ 
\les&\|u\|_{L^\infty}^2\Big(\|\nabla d\|_{L^2}^2+\|\nabla^2 d\|_{L^2}^2+\|\nabla^3 d\|_{L^2}^2\Big)
+\|\nabla d\|_{L^\infty}^2\Big(\|\nabla u\|_{L^2}^2+\|\nabla^2 u\|_{L^2}^2\Big)\\
&+\|\nabla u\|_{H^1}^2\|\nabla^2 d\|_{H^1}^2\le C.
\end{split}
\eeq
Similarly, for $L_3$, we have
\beq\label{nabla^4d-3}
\begin{split}
L_3\les& \||\nabla d|^2\|_{L^2}^2+\||\nabla d||\nabla^2d|\|_{L^2}^2
+\||\nabla^2 d|^2\|_{L^2}^2+\||\nabla d||\nabla^3d|\|_{L^2}^2\\ 
\les&\|\nabla d\|_{L^\infty}^2\|\nabla d\|_{H^2}^2+\|\nabla^2 d\|_{H^1}^4\le C.
\end{split}
\eeq
Substituting (\ref{nabla^4d-1})-(\ref{nabla^4d-3}) into (\ref{nabla^4d}), we have
\beq\label{nabla^4d-total}\begin{split}
\|\nabla^4d\|_{L^2}^2\les&\|\nabla^2 d_t\|_{L^2}^2+1.
\end{split}
\eeq
Integrating (\ref{nabla^4d-total}) over $(0,t)$, and using (\ref{nabla^2d_t-L^2}), we 
establish Corollary \ref{clccoro2.7}. \endpf\\

{\noindent\bf Proof of Theorem \ref{umaintheorem}:}

By the above estimates, we know that both (\ref{uniform_est1}) and (\ref{uniform_est2}) are valid. 
Hence $T_*$ is not the maximum
 time for the strong solution $(\rho, u, d)$.  This contradicts the
definition of $T_*$. The proof of Theorem
\ref{umaintheorem} is complete. \qed

\end{document}